\documentclass{article}
\usepackage[utf8]{inputenc}
\usepackage{graphicx,xcolor}
\usepackage{amsmath,amsfonts,amsthm,amssymb}
 \usepackage{hyperref}
\usepackage[english]{babel}
\renewcommand\theenumi{{(\roman{enumi})}}

\numberwithin{equation}{section}

\usepackage[a4paper,margin=1in]{geometry}
\newcommand{\borelset}{A}

\newtheorem{theorem}{Theorem}
\numberwithin{theorem}{section}
\newtheorem{thm}[theorem]{Theorem}
\newtheorem{prop}[theorem]{Proposition}
\newtheorem{lem}[theorem]{Lemma}

\newcommand{\LM}[1]{\hbox{\vrule width.2pt \vbox to#1pt{\vfill \hrule width#1pt height.2pt}}}
\newcommand{\LL}{{\mathchoice{\,\LM7\,}{\,\LM7\,}{\,\LM5\,}{\,\LM{3.35}\,}}}
\newcommand\calL{\mathcal{L}}
\newcommand\calH{\mathcal{H}}
\newcommand\sym{\mathrm{sym}}
\newcommand\qc{\mathrm{qc}}
\newcommand\rel{\mathrm{rel}}
\newcommand\Eel{E_\mathrm{el}}
\newcommand\Esur{E_\mathrm{sur}}
\newcommand\Aconstr{\mathcal A_\mathrm{constr}}
\newcommand\Aone{\mathcal A_{TB}}
\newcommand\Athree{\mathcal A_{LR}^\alpha}
\newcommand\ve{\varepsilon}
\newcommand\eps{\varepsilon}
\renewcommand\d{\partial}

\newcommand\Lip{\text{\rm Lip}}
\newcommand\Z{\mathbb{Z}}
\newcommand\R{\mathbb{R}}
\begin{document}

\begin{center}
{\Large
{\LARGE An energy minimization approach\\
to twinning with variable volume fraction}
}\\[5mm]
{\today}\\[5mm]
Sergio Conti$^{1}$, Robert V. Kohn$^{2}$ and Oleksandr Misiats$^{3}$\\[2mm]
{\em $^{1}$ 
 Institut f\"ur Angewandte Mathematik,
Universit\"at Bonn,\\ 53115 Bonn, Germany}\\[1mm]
{\em $^{2}$ Courant Institute, New York University, New York, NY, 10012}\\[1mm]
{\em $^{3}$ Department of Mathematics and Applied Mathematics, Virginia Commonwealth University,\\
Richmond, VA, 23284}\\[3mm]
    \begin{minipage}[c]{0.9\textwidth}
In materials that undergo martensitic phase transformation, macroscopic loading often leads to the creation and/or rearrangement of elastic domains. This paper considers an example {involving} a single-crystal slab made from two martensite variants. When the slab is made to bend, the two variants form a characteristic microstructure that we like to call ``twinning with variable volume fraction.''
Two 1996 papers by Chopra et. al. explored this example using bars made from InTl, providing considerable detail about the microstructures they observed. Here we offer an energy-minimization-based model that is motivated by their account.  It uses geometrically linear elasticity, and treats the phase boundaries as sharp interfaces. For simplicity, rather than model the experimental forces and boundary conditions exactly, we consider certain Dirichlet or Neumann boundary conditions whose effect is to require bending. This leads to certain nonlinear (and nonconvex) variational problems that represent the minimization of elastic plus surface energy (and the work done by the load, in the case of a Neumann boundary condition). Our results identify how the minimum value of each variational problem scales with respect to the surface energy density. The results are established by proving upper and lower bounds that scale the same way. The upper bounds are ansatz-based, providing full details about some (nearly) optimal microstructures. The lower  bounds are ansatz-free, so they explain why no other arrangement of the two phases could be significantly better.
\end{minipage}
\end{center}

\section{Introduction}

In materials that undergo martensitic phase transformation, large-scale
elastic deformation is typically accommodated by the creation or rearrangement
of small-scale elastic domains. This paper considers an example involving a
single-crystal slab made from two martensite variants. When the slab is made to
bend, the two variants form a characteristic microstructure that we like to
call ``twinning with variable volume fraction.'' Two 1996 papers \cite{ChoWut,Chopra-Roytburd-Wuttig}
explored this example using bars made from InTl, providing
considerable detail about the microstructures they observed. Here we
offer an energy-minimization-based model that is motivated by their account.
Our main accomplishments can be summarized as follows:

\begin{enumerate}
\item[(a)] We formulate an energy-minimization-based mathematical model capturing
many features of the scenario envisioned in \cite{ChoWut,Chopra-Roytburd-Wuttig} (which is summarized by
Figure \ref{fig:from-met-trans}). Actually, we formulate two such models: one in which the slab is
made to bend by specifying a displacement-type boundary condition at two opposite sides, and another in
which the slab is made to bend by applying a suitable boundary load. In each model,
the surface energy of the twin boundaries is a small parameter $\ve$. While our results
are valid for any $\ve \in (0,1]$, it is natural to focus on the limiting behavior as
$\ve \rightarrow 0$ since this corresponds to the presence of fine-scale
twinning as seen in the experiments.

\item[(b)] For each model, we provide an ansatz-based upper bound. It provides a candidate
domain structure, and a clear indication of the elastic strain associated with this
structure.

\item[(c)] For each model, we provide an ansatz-free lower bound, whose energy scaling law
(as $\ve \rightarrow 0$) matches that of our upper bound (modulo prefactors).
The lower bounds show rigorously that no domain structure can do significantly better
than those associated with our upper bounds. Moreover, the proofs of the lower bounds
provide intuition about why, in the presence of bending, the material forms the microstructure
that is seen.
\end{enumerate}

While the crystallography implicit in our elastic energy is motivated by the discussion
in \cite{ChoWut,Chopra-Roytburd-Wuttig}, we do not attempt to model the experimental forces or boundary
conditions. Rather, we consider certain Dirichlet or Neumann boundary conditions
whose effect is to induce bending similar to what is seen in the experiments.

Since we are interested in the limiting behavior as the surface energy becomes negligible,
it is natural to consider the \emph{relaxation} of the elastic energy (see Section \ref{sec:getting-started} for
its definition). The deformations with relaxed energy zero are those achievable with
arbitrarily small elastic energy, provided one places no restriction on the length scale of
the microstructure. These are, in a sense, the deformations achievable by mixing the two
variants. A particular minimizer $u^*$ of the relaxed energy (defined by \eqref{relaxed})
lies at the heart of our analysis. It captures our vision of what is happening in
the experiments, macroscopically. When we impose Dirichlet boundary conditions, they require
the some components of the deformation to agree with a multiple of $u^*$ on two opposite
faces of the domain. Our lower bounds rely on the fact that the relaxed solution $u^*$ is
achieved \emph{only} in the limit as the microstructural length scale tends to $0$;
including surface energy as well as elastic energy prevents infinitesimal phase mixtures,
and therefore increases the elastic energy.

Our ansatz-based upper bounds are rather elementary. They nevertheless carry useful
information, since they suggest specific microstructures that approach the behavior
of the relaxed solution. Our microstructures are presumably not optimal -- for example,
they do not solve any Euler-Lagrange equations. However our matching lower bounds
indicate that they approximately optimize the sum of elastic plus surface energy.

In our model with a Neumann boundary condition, it is natural to ask what the force-displacement
curve would look like. While our rigorous results address only the energy scaling law, our
ansatz-based upper bound suggests the existence of a critical load; below this load,
energy minimization prefers a single-domain state, while above this load, energy minimization
prefers the maximum bending that can be achieved using twinning with variable volume fraction
(see the discussion in Section \ref{subsec:main-results} just after Theorem \ref{theolowerboundneumann}).

To be clear: the goal of this paper is \emph{not} to provide a complete model of the
experiments reported in \cite{ChoWut,Chopra-Roytburd-Wuttig}. Rather, it is to begin the mathematial analysis of how
bending leads to twinning with variable volume fraction. Indeed, while twinning with constant
volume fraction has received a lot of attention, the analysis of problems involving
twinning with variable volume fraction has barely begun. (See Section \ref{subsec:related-work} for further
comments on the scientific context of this work.)

Another caveat: there is room for doubt about whether the experiments in \cite{ChoWut,Chopra-Roytburd-Wuttig} actually
involve mixtures of just two variants. Indeed, the later of the two papers (reference \cite{Chopra-Roytburd-Wuttig})
reports observing a microstructure with two distinct length scales (a mixture of ``polydomain
phases'' in the language of \cite{Roytburd-Wuttig-Zhukovskiy}, resembling what is known as a second-rank laminate in
the mathematical literature). Similar observations have been seen in recent work on the bending of
a cylinder made from NiMnGa \cite{Chu19}. We do not think such microstructures can be made by
mixing just two variants. It remains a challenge to model them, and to explain how
they are produced by the macroscopic bending of a slab or cylinder.

Despite the uncertain relationship between our two-variant model and the experiments, we believe the model is still
worthy of study. Indeed, its analysis provides fresh intuition about twinning with variable volume fraction,
and tools that will be useful in other settings (including perhaps models involving more than two variants).

This paper uses geometrically linear elasticity. It is natural to ask whether something
similar could be done in a geometrically nonlinear framework. The answer is that the constructions behind our upper
bounds have geometrically nonlinear analogues; however, we do not know how to prove corresponding lower bounds in
a geometrically nonlinear setting. As an initial step toward that goal, we recently considered a geometrically
nonlinear problem involving the bending of a two-dimensional bar made from two variants with a single
rank-one connection \cite{ckm-mmmas}.

Problems involving a scalar-valued unknown are often simpler than vector-valued analogues. It is therefore natural
to ask whether there is a scalar-valued analogue of twinning with variable volume fraction. A problem of that type
was considered briefly in \cite{KohnOtto97}, and it is studied more deeply in the forthcoming paper \cite{KohMulMis}. However,
the analogy between that scalar problem and what we do here is far from
complete. Indeed, in the present setting (with our Dirichlet boundary conditions) the relaxed problem has minimum energy
zero, while in the scalar setting of \cite{KohnOtto97,KohMulMis} the minimum of the relaxed energy is strictly positive.

The rest of this paper is organized as follows: Section \ref{sec:getting-started} presents our model, states our main results, and
puts our work in context by discussing some related literature. Section \ref{sec:upper-bounds}, which is relatively short
and elementary, gives the proofs of our upper bounds. Section \ref{sec:lower-bounds}, which occupies roughly half the
paper, gives the proofs of our lower bounds; a guide to its structure will be found at the
beginning of that section.

\begin{figure}
 \begin{center}
Figures not included for copyright reasons. See Figure 3b and Figure 4 of \cite{Chopra-Roytburd-Wuttig}
\end{center}
\caption{Left: Micrograph of microstructure in a sample of InTl shape-memory alloy under bending, from
\cite[Fig.~3b]{Chopra-Roytburd-Wuttig}. Two variants are apparent; the volume fraction has a strong dependence on the vertical coordinate.
Right: Geometry of the experiment and schematic representation of the observed microstructure from \cite[Fig.~4]{Chopra-Roytburd-Wuttig}.
The micrograph corresponds to the sketch g.
The corresponding microstructure in our model is presented in Figure~\ref{fig:UB} below.
{Reprinted by permission from Springer Nature from
Metallurgical and Materials Transactions A, \textbf{27}, 1695--1700,
  (1996), Copyright 1996, Ref.~\cite{Chopra-Roytburd-Wuttig}.}}
\label{fig:from-met-trans}
\end{figure}

\section{Getting started} \label{sec:getting-started}

\subsection{Problem set up} \label{subsec:problem-setup}

The right hand side of Figure \ref{fig:from-met-trans} provides a sketch of the twinning with varying volume
fraction that the authors of \cite{ChoWut,Chopra-Roytburd-Wuttig} reported as the result of bending an
InTl slab with well-chosen crystallography. As one sees from parts (e)--(h) of the sketch, in the bent
slab, one variant predominates on the top while the other predominates on the bottom. This is achieved by
the boundaries between the two variants tilting slightly from their stress-free orientations in the
unbent slab, shown in parts (a)--(d) of the figure. Note that the sketch is not to scale; the length
scale of the twinning is actually quite small, and the angle of the tilt is correspondingly small,
as one sees from the left hand side of Figure \ref{fig:from-met-trans}. (We refer to the
sample as a slab, though in fact it was a cylinder with a rectangular cross-section. One sample was
$1.735 {\rm cm}$ long with a $.111 {\rm cm} \times .133 {\rm cm}$ cross-section, and the other
was similar.)

The paper \cite{ChoWut} takes the view that this microstructure is a mixture of two martensite variants. Our goal is {to}
formulate a quantitative model based on that idea, and to explore its consequences. Since the variants come from
a cubic-tetragonal phase transformation and we are using geometrically linear elasticity, it would be standard
(see e.g. \cite{Bhat}) to assume that their stress-free strains are
\[
 \pm
\left(
  \begin{array}{ccc}
    0 & 0 & 0 \\
    0 & 1 & 0 \\
    0 & 0 & -1 \\
  \end{array}
\right) .
\]
With this choice the twin planes would be normal to $(0,1,1)$ and $(0,1,-1)$. The sketch in
Figure \ref{fig:from-met-trans} uses this choice (indeed, part (a) indicates that the twin
plane is normal to $(0,1,-1)$ and {the} slab is longest along the axis $(0,0,1)$). For us, however, it is
convenient to work in a rotated coordinate system where the stress-free strains are
\begin{equation}\label{wells}
\pm \left(
  \begin{array}{ccc}
    0 & 0 & 0 \\
    0 & 0 & 1 \\
    0 & 1 & 0 \\
  \end{array}
\right)
\end{equation}
and the twin planes are normal to $(0,1,0)$ and $(0,0,1)$. This coordinate system is related to the
original one by a rotation with angle $\pi/4$ about the $(1,0,0)$ axis. Consistent with the sketch, we
shall assume that in absence of bending the interfaces between the two variants are normal
to $(0,0,1)$ in our rotated coordinate system. This choice of coordinates is convenient because
it permits our upper-bound ansatz to be invariant in the $x_2$ direction.
\medskip

\noindent {\sc Our spatial domain.} By elasticity scaling (replacing $u(x)$ by $\lambda u(x/\lambda)${)},
only the shape of the domain matters, not its actual dimensions. Throughout this paper, our spatial
domain will be the unit cube
$$
\Omega:=(-1,1)^3 .
$$
Our upper bounds involve microstructures that are periodic in $x_3$ and invariant in $x_2$,
so they are easily extended to slab-like domains like those sketched in
Figure \ref{fig:from-met-trans}. Our lower bounds surely extend to rectangular solids that
are not cubes, however for domains that are highly eccentric (slabs or cylinders) it would take
extra work to identify the bounds' dependence on length vs height vs width, and we are not sure
the result would have the same scaling as the upper bounds when the domain is highly eccentric.
\bigskip

\noindent {\sc Our elastic energy.} For any $u\in W^{1,2}(\Omega;\R^3)$ we define
\begin{equation}\label{strain}
 e(u):=\frac12 (Du+Du^T)\,, \quad
e_{ij}(u)= \frac{\d_{x_j} u_i + \d_{x_i} u_j}{2}, \hskip2mm i,j = 1,2,3.
\end{equation}
In view of (\ref{wells}), we impose the constraint that
\begin{equation}\label{constraint}
e_{23}(u) = \pm 1 \hskip5mm\text{ almost everywhere}.
\end{equation}
Our elastic energy enforces this constraint: it is
\begin{equation} \label{elastic-energy}
\Eel[{u}] : = \int_{\Omega}  W(e(u)) d{x},
\end{equation}
where the elastic energy density $W:\R^{3\times 3}_{\sym}\to[0,\infty]$ is defined by
\begin{equation} \label{elastic-energy-density}
W(\xi):=
  \begin{cases}
  \xi_{11}^2+\xi_{22}^2+\xi_{33}^2+2\xi_{12}^2+2\xi_{13}^2, & \text{ if }  |\xi_{23}|=1,\\
  \infty,&\text{ otherwise.}
  \end{cases}
\end{equation}

Some readers might wonder why we impose $e_{23} = \pm 1$ as a constraint, rather than taking the
elastic energy density to be the minimum of two quadratic functions, for example
\begin{equation}\label{elastic-energy-quadratic}
\min
\left\{
\left\| e(u) - \left(\begin{smallmatrix} 0 & 0 & 0\\ 0 & 0 & 1 \\ 0 & 1 & 0 \end{smallmatrix}\right) \right\|^2 ,
\left\| e(u) + \left(\begin{smallmatrix} 0 & 0 & 0\\ 0 & 0 & 1 \\ 0 & 1 & 0 \end{smallmatrix}\right) \right\|^2
\right\}.
\end{equation}
The answer is that the choice \eqref{elastic-energy-density} is mathematically convenient, since it permits us
to identify the interface between the two variants as the surface where $e_{23}$ changes sign. We believe that
use of \eqref{elastic-energy-quadratic} in place of \eqref{elastic-energy} would not fundamentally change the
results. Elastic energies analogous to our $W$ have been used in many studies of how elastic and surface
energy interact to set the geometry and length scale of twinning in martensitic phase transformation (an early
example is \cite{KohMul94}, and a more recent example with many references is \cite{ConDierMelchZwick}).
\bigskip

\noindent {\sc Our surface energy.}
Recalling \eqref{constraint}, we define the surface energy as
\begin{equation}\label{eqdefesur}
\Esur[u]:= \int_{\Omega} |D (e_{23}(u))| d {x},
\end{equation}
which is twice the measure of the interface between $e_{23} = 1$ and $e_{23} = -1$.
As usual in this context, the integral in
\eqref{eqdefesur} should be interpreted as the total variation of the measure
$D (e_{23}(u))$
if $e_{23}(u)\in BV(\Omega)$, and $\infty$ otherwise.

Since the experiments we want to model involve fine-scale twinning, we expect the surface energy density to
be small. Thus the energy of our sample is the sum of elastic energy and a small parameter $\ve >0$ times
surface energy:
\begin{equation}\label{energy3D}
E_\ve[{u}]:= \Eel[{u}] + \ve \Esur[{u}].
\end{equation}
While our results are valid for any $\ve \in (0,1]$, they are mainly of interest in the limit $\ve \rightarrow 0$.
\bigskip

\noindent {\sc The relaxed energy, and our relaxed solution.}
It is easy to see that the quasiconvex envelope of our elastic energy
density $W$ is the convex function
\begin{equation}
  W^\qc(\xi):=
  \begin{cases}
  \xi_{11}^2+\xi_{22}^2+\xi_{33}^2+2\xi_{12}^2+2\xi_{13}^2, & \text{ if }  |\xi_{23}|\le 1,\\
  \infty,&\text{ otherwise.}
  \end{cases}
 \end{equation}
(The proof uses lamination, taking advantage of the fact that $e_2\otimes e_3$ is rank-one.)
In particular, the relaxation of $\Eel$ is the functional
\begin{equation}
\Eel^\rel[{u}]  = \int_{\Omega}  W^\qc(e(u)) d{x}.
\end{equation}
Its minimizers are the weak limits of minimizing sequences of \eqref{elastic-energy}. In particular, a deformation
with relaxed energy $0$ is one that can be achieved by fine-scale twinning with arbitrarily small elastic energy.

The following minimizer of the relaxed problem plays a crucial role in our analysis:
\begin{equation}\label{relaxed}
    u^*(x):=
    \begin{pmatrix} - x_2 x_3  \\
       x_1 x_3 \\
      x_1 x_2
    \end{pmatrix}.
  \end{equation}
One verifies by elementary calculation that the relaxed energy vanishes at $u = \alpha u^*$ whenever $|\alpha| \leq 1$.
\emph{The central thesis of this paper is that in the limit $\ve \rightarrow 0$, the bending seen in the experiments
corresponds macroscopically to a deformation of this form.}

This deformation deforms the slab shown in Figure \ref{fig:from-met-trans} into a saddle shape. Indeed, the midplane of
the slab is the $x_2,x_3$ plane. Taking $\alpha =1$ for simplicity, the normal displacement of this plane is
$u_1^* (x_1,x_2) = - x_2 x_3 = \frac{1}{4}(x_2 - x_3)^2 - \frac{1}{4} (x_2 + x_3)^2$. Remembering that we are working
in a rotated coordinate system relative to that of Figure \ref{fig:from-met-trans}, the principal axes of this bending
deformation are precisely those parallel and perpendicular to the long axes of the sample.

The experimental papers discuss ``bending'' but do not mention seeing such a saddle shape. This is, we presume,
because their samples were much longer than they were wide -- really, more like cylinders than slabs. Therefore the
bending in the long direction would have been pronounced, while the opposite-oriented bending in the {orthogonal} direction
would hardly have been noticeable. The prediction of a saddle shape makes physical sense. Indeed, our model
assumes that the slab is made from two variants of martensite, which come from a volume-preserving
(cubic-tetragonal) phase transformation; therefore, after deformation the volume of the upper half
of the slab should match that of the lower half. This requires the image of the slab's midplane to have
mean curvature zero.
\bigskip

\noindent {\sc Our Dirichlet boundary conditions.}
As mentioned already in the Introduction, we do not attempt to directly model the experiments; rather, we shall
make our sample bend by imposing suitable Dirichlet or Neumann boundary conditions. We discuss the Dirichlet conditions
here; the Neumann conditions will be introduced a little later, in Section \ref{subsec:main-results}.

Our Dirichlet conditions require that certain components of $u$ agree with a relaxed solution on two
opposite faces. We consider two distinct alternatives. The first imposes conditions on the top and
bottom faces $x_1 = \pm 1$:
\begin{equation}\label{BC alpha=1}
    \begin{cases}
    u_2(\pm 1, x_2, x_3) = \pm x_3,\\
    u_3(\pm 1, x_2, x_3) = \pm x_2 .
    \end{cases}
  \end{equation}
(No condition is imposed on $u_1$.) This amounts to requiring that $u_2$ and $u_3$ agree with the relaxed solution
$u^*$ at $x_1 = \pm 1$. It is crucial here that we use $u^*$, not $\alpha u^* $ for some $|\alpha| < 1$. Indeed,
since $e_{23}(u^*) =1$ when $x_1 = 1$ and $e_{23}(u^*) = -1$ when $x_1 = -1$, $u^*$ describes a bent
configuration whose top and bottom faces are not infinitesimally twinned -- rather, each consists of a pure variant.
(The sketch in Figure \ref{fig:from-met-trans} corresponds to $\alpha$ close to but not quite equal to $1$.)
The corresponding class of test functions will be denoted by $\mathcal{A}_{TB}$,
\begin{equation}\label{class1}
\mathcal{A}_{TB} := \{u\in W^{1,2}(\Omega;\R^3): u \text{ satisfies  (\ref{BC alpha=1})} \}.
\end{equation}

It is not so natural to impose $u = \alpha u^*$ on the top and bottom boundaries for $|\alpha|<1$, since
this would require the top and bottom faces to twin infinitesimally (whereas when surface energy is included we
expect the length scale of twinning to be strictly positive).

Our second alternative specifies the value of $u_3$ on the left and right faces $x_3 = \pm 1$.
It has the advantage of working equally well for any $\alpha \in [-1,1]$; in other words, we can
impose any amount of bending rather than just the maximal amount. The condition we impose under this
alternative is
\begin{equation}\label{BC alpha<1}
    u_3(x_1, x_2, \pm 1) = \alpha x_1 x_2 \, .
\end{equation}
The corresponding class of test functions will be denoted by $\mathcal{A}_{LR}^\alpha$:
\begin{equation}\label{class3}
\mathcal{A}_{LR}^\alpha := \{u\in W^{1,2}(\Omega;\R^3): u \text{ satisfies (\ref{BC alpha<1})} \}.
\end{equation}

Let us dwell a bit on the value of permitting $|\alpha| < 1$. The relaxed solution $u = \alpha u^*(x)$ has
$$
e(u) = \begin{pmatrix} 0 & 0 & 0\\ 0 & 0 & \alpha x_1\\ 0 & \alpha x_1 & 0 \end{pmatrix}.
$$
Since $\alpha x_1 = \theta \cdot 1 + (1-\theta) \cdot (-1)$ when $\theta = (1 + \alpha x_1)/2$, this deformation
corresponds microscopically to a mixture of the two variants with a volume fraction that varies linearly in
$x_1$ (as shown in Figure \ref{fig:from-met-trans}, and as will be clear from the constructions that prove our
upper bounds). At the bottom ($x_1 = -1${)} we expect volume fraction $(1-\alpha)/2$ of the $e_{23} = +1$ phase and
volume fraction $(1+\alpha)/2$ of the $e_{23} = -1$ phase; at the top the volume fractions are reversed. The
macroscopic bending is obviously proportional to $\alpha$. In particular, the microstructure of the unbent slab
(shown in part (a) of the sketch in Figure \ref{fig:from-met-trans}) corresponds to $\alpha = 0$.

We remark that for any $\ve > 0$ and any $\alpha\in[-1,1]$, our elastic plus surface energy functional $E_\eps$
(defined by \eqref{energy3D}) achieves its minimum over either $\mathcal{A}_{LR}^\alpha $ or $\mathcal{A}_{TB}$.
Indeed, for any minimizing sequence $u_j$ the fact that $e_{23}(u_j)$ is bounded in $BV(\Omega)$
implies that it has a strongly converging subsequence, hence the constraint $e_{23} = \pm 1$ passes to the
limit. The dependence on the other terms is convex, and the boundary data pass to the limit by the trace theorem
(since a bound on the elastic energy implies a bound on the $W^{1,2}$ norm of $u$).

We close this discussion by introducing one more class of test functions. In Section \ref{sec:lower-bounds}, where we present our
lower bounds, we will offer two distinct approaches to the lower bound for $u \in \mathcal{A}_{TB}$. One of them
uses a duality argument and assumes the additional structural hypothesis that
\begin{equation}\label{ansatz}
    \begin{cases}
      u_1 = u_1(x_1,x_2,x_3); \\
      u_2 = u_2(x_1,x_3); \\
      u_3 = x_2 \hat{u}_3(x_1, x_3)
    \end{cases}
  \end{equation}
for some function $\hat{u}_3$. (This is consistent with our upper bound constructions, and with the sketch in Figure
\ref{fig:from-met-trans}.) The corresponding set of deformations will be called $\Aconstr$:
\begin{equation}\label{classc}
\Aconstr:= \{u\in W^{1,2}(\Omega;\R^3): {u} \text{ satisfies
\eqref{ansatz}} \}.
\end{equation}
Our second proof of the lower bound for $u \in \mathcal{A}_{TB}$ has a broader scope than the duality-based argument,
since it doesn't require the structural hypothesis \eqref{ansatz}. However our two proofs of the lower bound
provide somewhat different insight about how the inclusion of surface energy affects the length scale and character
of the microstructure. Therefore it seems useful to provide them both.

\subsection{Main results} \label{subsec:main-results}
We are now in position to formulate our main results. The uniqueness of the relaxed solution (subject to either of
our Dirichlet boundary conditions) lies, as already noted, at the heart of our analysis:

\begin{prop}\label{P1:rigid}
The function $u^*$, defined by \eqref{relaxed}, has the following properties:
\begin{enumerate}
\item\label{P1:rigidone} The function $u^*$ minimizes the relaxed problem
\[
\min\{\Eel^\rel[{u}] : u\in \Aone\}.
\]
Any other minimizer has the form $u^*(x)+ce_1$ for some $c\in\R$.

\item\label{P1:rigidthree} For any $\alpha\in[-1,1]$, the function $\alpha u^*$ minimizes the problem
\[
\min\{\Eel^\rel[{u}] : u\in \Athree \}.
\]
Any other minimizer has the form
\begin{equation}\label{eqriglr}
u(x)=\alpha u^*(x)+\begin{pmatrix}
                      dx_2+c\\
                      -dx_1+\psi(x_3)\\
                      0
                     \end{pmatrix}
\end{equation}
for some $c,d\in\R$ and $\psi\in\Lip((-1,1))$ with $\Lip(\psi)\le 2(1-|\alpha|)$.
\end{enumerate}
\end{prop}
The proof is given in Section~\ref{subsec:LB-1} below. We shall discuss in a moment why this result is important.

We now state our main results concerning the energy scaling laws using Dirichlet boundary conditions. When the
constraints are at the top and bottom boundaries, our main result is as follows.

\begin{thm}\label{thm1}
There are $c_1,c_2>0$ such that for any $\eps\in(0,1]$
\begin{equation}
c_1 \ve^{2/3} \leq
\min \{ E_\ve[{u}]: u\in \mathcal{A}_{TB}\}
\leq c_2 \ve^{2/3}.
\end{equation}
\end{thm}
\begin{proof}
The upper bound follows from
Proposition~\ref{propubuniversal} with $\alpha=1$, and
the lower bound from Proposition~\ref{proplowerbound}\ref{theolowerboundbcx1}.
\end{proof}

When the constraints are at the left and right boundaries, our main result is as follows.

\begin{thm}\label{thm2}
There are
$d_1, d_2>0$ such that for any $\alpha\in[-1,1]$ and any $\eps\in(0,1]$
\begin{equation}
d_1 \min\{|\alpha|^2,|\alpha|^{2/3}\ve^{2/3}\} \leq
\min \{ E_\ve[{u}]: u\in \mathcal{A}_{LR}^\alpha\} \leq d_2 \min\{|\alpha|^2,|\alpha|^{2/3}\ve^{2/3}\}.
\end{equation}
\end{thm}
\begin{proof}
The upper bound follows from
Proposition~\ref{propubuniversal}, and
the lower bound from  Proposition~\ref{proplowerbound}\ref{theolowerboundbcx4}.
\end{proof}

For each of the preceding theorems, the upper bound uses an explicit construction reminiscent of those
experimentally observed. This construction is described in Section \ref{sec:upper-bounds}. The lower bounds draw inspiration
from Proposition \ref{P1:rigid}, which shows in essence that to drive the elastic energy to $0$ one needs
twinning on an infinitesimal length scale. When $\ve > 0$ the surface energy prevents this. Our task in proving the
lower bounds is to assess the excess elastic energy that must be present due to the surface-energy-imposed
limitation on the length scale of twinning. One argument, presented in Section \ref{subsec:LB-2}, relies on the
convexity of the relaxed energy. A second argument, presented in Sections \ref{subsec:LB-3}--\ref{subsec:LB-5},
relies on a quantitative
version of Proposition \ref{P1:rigid} -- specifically, Proposition \ref{proplowerstructure} -- which shows
that if a map $u$ has small relaxed energy than it is in a certain sense close to $\alpha u^*$.

We turn now to our results using a Neumann boundary condition. The idea is that one can induce bending by
imposing suitable forces at the left and right hand boundaries $x_3 = \pm 1$. Precisely, we consider for
some $\gamma\in\R$ a force of the type $f'(x_1,x_2,\pm1):=\pm \gamma x_2 e_1$ acting on the faces normal
to $e_3$. This type of force can be represented by the linear functional
$M':W^{1,2}(\Omega;\R^3)\to\R$ defined by
\begin{equation}\label{eqdefM1}
M'(u):=\int_{(-1,1)^2}
x_2 (u_1(x_1,x_2,1)-u_1(x_1,x_2,-1)) d\calL^2(x_1,x_2).
\end{equation}
The Neumann problem then amounts to minimizing $E_\eps[u]-\gamma M'(u)$ over all $u\in W^{1,2}(\Omega;\R^3)$.
Our main result in this case is the following.

\begin{thm}\label{theolowerboundneumann}
There is {$c>0$} such that for all $\eps\in (0,1]$ and all $\gamma\in\R$
\begin{multline*}
c\min\{-\gamma^2, \frac1{c^2}\eps^{2/3}-|\gamma|\}\le \inf\{E_\eps[u]-\gamma M'(u): u\in W^{1,2}(\Omega{;\R^3})\} 
  \le \frac1c\min\{-\gamma^2,  c^2\eps^{2/3}-|\gamma|\}.
\end{multline*}
\end{thm}

\begin{proof}
The upper bound follows from Proposition~\ref{propubneu}, and the lower bound from
Proposition~\ref{proplbneu}. We remark that the constant $c$ is common. Indeed, if $0<c_1\le c_1'$,
$0<c_2\le c_2'$, $0<c_3\le c_3'$, then for all $\gamma\in\R$, $\eps>0$ one has
\begin{equation}\label{eqneumonot}
\min\{-\frac1{c_1 }\gamma^2, c_2\eps^{2/3}-\frac1{c_3}|\gamma|\}\le
\min\{-\frac1{c_1'}\gamma^2, c_2'\eps^{2/3}-\frac1{c_3'}|\gamma|\}.
\end{equation}
\end{proof}

\begin{figure}
\begin{center}
\includegraphics[width=4cm]{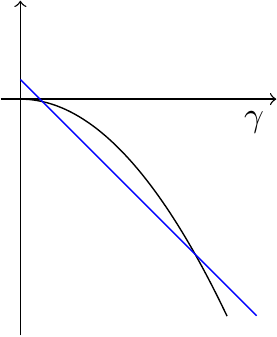}
\end{center}
\caption{Sketch of the function which estimates  the optimal energy in
Theorem~\ref{theolowerboundneumann}. The parabola is the elastic solution, the
straight line the one with microstructure.}
\label{figneumann}
\end{figure}

The dependence of the minimum energy on $\gamma$ is interesting. As illustrated in
Figure~\ref{figneumann}, there is a competition between two different response regimes. One involves
the elastic response of the single-variant state; it results in an energy of order $- \gamma^2$. The
other involves microstructure -- in fact, fully-developed bending of the type associated with relaxed energy
$\alpha u^*$ with $\alpha = 1$ when $\gamma < 0$ and $\alpha = -1$ when $\gamma > 0$ (see the proof of
Proposition \ref{propubneu}) -- and its energy is of order $c \ve^{2/3} - |\gamma|$. The regime with microstructure
has smaller energy when $\eps^{2/3}\lesssim |\gamma|\lesssim 1$. Thus, energy-minimization suggests that in
a load-controlled experiment, the sample would initially respond elastically, then suddenly bend sharply
when $\gamma$ exceeds a critical value that's of order $\ve^{2/3}$.

While we have restricted our attention to the load associated with the boundary term \eqref{eqdefM1},
we believe that similar results could be obtained (by similar arguments) for more general forces that have
the same symmetry, and for forces acting on other faces.

\subsection{Related work} \label{subsec:related-work}

As already explained in the Introduction, the goal of this paper is to begin the analysis of how bending induces
twinning with variable volume fraction, in a slab or cylinder made from a material with two symmetry-related
variants.

While our work is specifically motivated by \cite{ChoWut,Chopra-Roytburd-Wuttig}, the idea that bending is associated
with twinning with variable volume fraction is much older; see, for example,
\cite{Basinski-Christian,Otsuka,Roytburd-Wuttig-Zhukovskiy}.
While these papers offers some analysis based on energy minimization, they provide only upper bounds, and their
dislocation-based approach is very different from the one in this paper.

Phenomena similar to twinning with variable volume fraction arise in some situations that have nothing to do
with bending. In particular, ``zigzag domain boundaries'' are seen in some ferroelectric--ferroelastic materials, see e.g.
\cite{FliHaa,MeeAul,Roitburd-1988}. In this setting the volume fraction varies on a macroscopic length
scale due to a nonlocal effect other than bending. Some simulations using a time-dependent Ginzburg-Landau
model were presented in \cite{Kuroda-etal-1998}, however we are not aware of any analysis comparable to that of the
present paper.

Something similar can also be found in certain optimal design problems. Indeed, shape optimization
sometimes leads to laminated microstructures whose volume fractions can vary macroscopically; see for example
Section 1.3 of \cite{Jouve-2014} and Figures 2 and 8 of \cite{Collins-Bhattacharya}. While the asymptotic effect
of including surface energy has not been studied in such settings, it has been considered in some shape optimization
problems where the volume fractions are constant \cite{Kohn-Wirth-uniaxial,Kohn-Wirth:15}.

As already noted in the Introduction, it is not entirely clear that the two-variant picture developed here is
an adequate description of the experiments in \cite{ChoWut,Chopra-Roytburd-Wuttig}. Indeed, the earlier paper
\cite{ChoWut} provides an account based entirely on mixtures of two variants, and the present paper is
based on that account. However the later paper \cite{Chopra-Roytburd-Wuttig} reports observing a
microstructure with two distinct length scales, resembling what is known as a second-rank laminate in
the mathematical literature. It remains a challenge to explain why bending produced such a structure.
Microstructures with two distinct length scales have also been observed when bending cylinders
made from NiMnGa \cite{Chu19}.

A phenomenon that resembles twinning with variable volume fraction -- but is actually quite different -- is
studied in \cite{Ganor-etal}. When two elastic phases come from a martensitic phase transformation, they have
two possible twin planes. Using both twin planes, one can arrange the phases in a cylindrical reference region
in such a way that the martensitic transformation makes the cylinder bend. This phenomenon is different from that
of the present paper because (i) it requires using both twinning systems, (ii) it does not involve microstructure, and
(iii) both before and after the phase transformation the elastic energy is exactly zero.

The mathematical context of our work is the use of energy minimization to study how the microstructure of an elastic
material with two or more variants changes in response to deformation or loading. Without attempting a survey of work
in this area, we mention the foundational studies \cite{Kha83,Roi} using geometrically linear
elasticity and \cite{BJ87,Ball-James-proposed-tests} using a geometrically nonlinear approach.
The monograph \cite{Bhat} provides a good introduction.

The minimization of elastic energy alone often predicts the formation of infinitesimal mixtures, whereas in
real materials we see mixtures with a well-defined length scale (at least locally). It is widely accepted that
the inclusion of surface energy sets a length scale, and in some settings also selects a preferred microstructural
pattern. This phenomenon has been analyzed quite extensively in studies of how two variants of martensite twin in the
vicinity of an austenite interface; see \cite{KohMu-philmag,KohMul94} for early work of this type, and
\cite{ChaCon15,ConDierMelchZwick} for recent contributions with many references. The microstructures
studied in that work are quite different from the ones considered here; in particular, they do not involve
twinning with varying volume fraction.

Our main results are upper and lower bounds for the minimum energy as a function of the surface energy density $\ve$.
The upper bounds come from explicit test functions. The lower bounds are more subtle -- their proofs use rigidity
theorems and/or convexity of the relaxed energy -- but they \emph{make no use of the Euler-Lagrange equations} that
characterize critical points of our functional. It is, in fact, difficult to use the stationarity or minimality of
elastic plus surface energy; however minimality has been used successfully in \cite{Con,ConDierKoserZwick}.

%\color{black}

\section{Upper bounds} \label{sec:upper-bounds}

\begin{prop}\label{propubuniversal}
{There is $c>0$ such that}
for every $\eps\in(0,1]$ and $\alpha\in[-1,1]$ there is $u_{\eps,\alpha}\in \Athree\cap\Aconstr$ with
  \begin{equation}\label{UpBndEl}
  E_\ve[{u_{\eps,\alpha}}] \leq c \min\{\alpha^2,|\alpha|^{2/3}
  \ve^{2/3}\}.
  \end{equation}
   Furthermore, if $\alpha = 1$, then $u_{\eps,1}\in \Aone \cap\Aconstr$ {and $u_{\eps,1}=u^*$ on $\partial\Omega$}.
  \end{prop}
\begin{proof}
We first remark that $\alpha^2\le |\alpha|^{2/3} \ve^{2/3}$ 
is equivalent to $\alpha^2\le \eps$. We shall use two different constructions in the two regimes.
  
We first consider the case $\eps\le\alpha^2$. Fix $n \geq 1$ to be chosen below, and denote $\delta:=\frac{1}{n}$. Let us partition the cube ${\overline\Omega} = [-1,1]^3$ into $2n$ slices in the direction $x_3$:
\[
{\overline\Omega} = \bigcup_{k=-n}^{n-1} \Omega_k, \text{ where } \Omega_k := [-1,1] \times [-1,1] \times [k\delta, (k+1)\delta].
\]
The upper bound is proved by construction of a test function $u_{\eps,\alpha}$ of the form
  \begin{equation}\label{eqdefuepsalpha}
  u_{\eps,\alpha}(x):=\begin{pmatrix}
  -\alpha x_2x_3 \\ u_2(x_1,x_3)\\ \alpha x_1x_2
                 \end{pmatrix}.
  \end{equation}
The function $u_2(x_1,x_3)$ is chosen to satisfy the following conditions:
\begin{itemize}
\item[{[i]}] $u_2(x_1,x_3) = u_2^*(x_1, x_2, x_3) = \alpha x_1 x_3$ for $x_3 = k\delta$, {$k\in\Z\cap[-n,n]$};
\item[{[ii]}] $e_{23}(u) = \pm 1$.
\end{itemize}
{Let us present a first} construction of $u_2$ in one domain $\Omega_k$
{(this will later be used if} $k$ is even, {as discussed after \eqref{el_bnd}}). Assume also that locally near the $\{x_3 = k\delta\}$ boundary the condition [ii] with $e_{23}(u) = 1$ holds. In view of the choice $u_3 = \alpha x_1 x_2$, near the $\{x_3 = k\delta\}$ boundary we  define
\[
u_2^{L}(x_1,x_3):= (2-\alpha x_1) x_3 + c_k^L(x_1).
\]
The condition [i] at $x_3 = k\delta$ immediately allows us to determine
\[
c_k^L(x_1) = 2 \alpha x_1 k \delta - 2 k \delta,
\]
and thus
\[
u_2^{L}(x_1,x_3):= (2-\alpha x_1) x_3 + 2 \alpha x_1 k \delta - 2 k \delta.
\]
We remark that this construction is equivalent to setting $u_2^L(x_1,x_3)=u_2^*(x_1,\cdot,k\delta)+(2e_{23}(u)-\partial_2u_3)(x_3-k\delta)$, with $e_{23}(u)=1$.

Assume now that locally near the boundary $x_3 = (k+1)\delta$,  the condition [ii] with $e_{23}(u) = -1$ holds. In a similar way, using the condition [i] at $x_3 = (k+1) \delta$, we may define
\[
u_2^{R}(x_1,x_3):= (-2-\alpha x_1) x_3 + 2 \alpha x_1 (k+1) \delta + 2 (k+1) \delta.
\]
Setting $u_2^{L}(x_1,x_3) = u_2^{R}(x_1,x_3)$, we obtain the interface equation
\begin{equation}\label{interface}
{x_3 = f_k^\alpha(x_1), \hskip1cm f_k^\alpha(x_1):=\frac{\alpha x_1 +1}{2}\delta  + k\delta.}
\end{equation}
{Obviously $k\delta\le f_k^\alpha(x_1)\le (k+1)\delta$ for all $x_1\in[-1,1]$.}
Altogether, in $\Omega_k$ we define
\begin{equation}\label{def_of_u2}
u_2(x_1,x_3) :=
\begin{cases}
u_2^{L}(x_1,x_3), \text{ if } {k\delta\le x_3 \leq f_k^\alpha(x_1);}\\
u_2^{R}(x_1,x_3), \text{ if } { f_k^\alpha(x_1)< x_3\le (k+1)\delta}.
\end{cases}
\end{equation}
It is straightforward to verify that this test function is continuous in $\Omega_k$, satisfies [i], [ii], and $e_{11} \equiv e_{22} \equiv e_{33} \equiv e_{13} \equiv 0.$ {By the condition [i] it is continuous in $\Omega$, and by construction $e_{23}=\pm 1$ almost everywhere. Similarly, from the definition we obtain $u_{\eps,\alpha}\in\Athree\cap\Aconstr$.}

We next consider the boundary conditions for $\alpha=1$. First,
\eqref{interface} yields that
{$f_k^1(1)=(k+1)\delta$ (i.e. the right edge of $\Omega_k$)  and
$f_k^1(-1)=k\delta$}
(i.e. the left edge of $\Omega_k$),
as illustrated in {Figure \ref{fig:UBb}}.
{Further, one easily checks that $u_2(1,x_3)=u_2^L(1,x_3)=x_3$
and $u_2(-1,x_3)=u_2^R(-1,x_3)=-x_3$ for
$x_3\in [k\delta, (k+1)\delta]$. Therefore 
{$u_{\epsilon,1}=u_*$ on $\partial\Omega$ and}
$u_{\epsilon,1}$} belongs to $\Aone$.

\begin{figure}
    \centering
    \includegraphics[height=4cm]{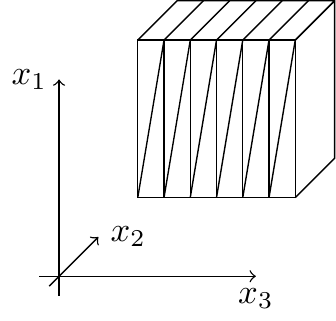}\hskip1cm
    \includegraphics[height=4cm]{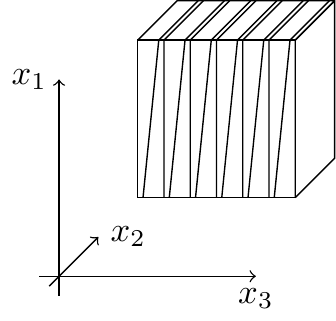}
    \caption{Sketch of the construction in the proof of Proposition~\ref{propubuniversal} for $\alpha = 1$ (left) and $\alpha\in(0,1)$ (right).}
    \label{fig:UBb}
\end{figure}
\begin{figure}
    \centering
    \includegraphics[height=4cm]{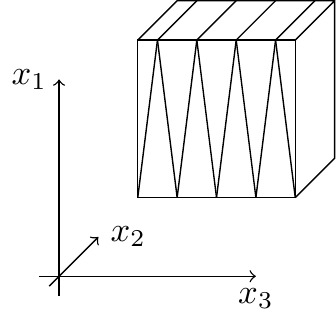}\hskip1cm
    \includegraphics[height=4cm]{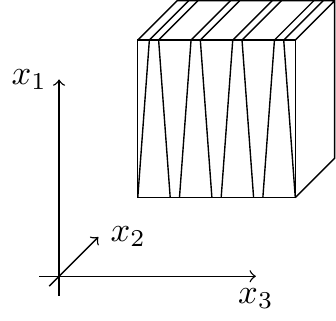}
    \caption{Sketch of the symmetric construction for $\alpha = 1$ (left) and $\alpha\in(0,1)$ (right).}
    \label{fig:UB}
\end{figure}

It remains to estimate $e_{12}$. Evaluating this strain separately for $u_2^{L}$ and $u_2^{R}$, we have
\[
e_{12}^L := \frac{1}{2}(\partial_2{u_1} + \partial_1{u_2^L}) = \alpha(k \delta - x_3),
\]
and
\[
e_{12}^R := \frac{1}{2}(\partial_2{u_1} + \partial_1{u_2^R}) = \alpha((k+1) \delta - x_3).
\]
Since $k\delta \leq x_3 \leq (k+1)\delta$, in both cases
\begin{equation}
|e_{12}(u_{\eps,\alpha})| \leq \alpha \delta \text{ a.e. in $\Omega_k$},
\end{equation}
{and therefore the same holds in $\Omega$. This leads to an estimate of the elastic energy,}
{\begin{equation}\label{el_bnd}
{\Eel[u_{\eps,\alpha}]} \le 8\alpha^2  \delta^2=\frac{8\alpha^2}{n^2}.
\end{equation}

{If this procedure is 
used for all values of $k$ (even and odd) one obtains the valid construction}
 depicted in Figure~\ref{fig:UBb}. However,
a more symmetric construction can be obtained
{using the above procedure for even $k$ and a}
 symmetric reflection for odd $k$,
which avoids the vertical interface for $e_{23}$ at $x_3 =k$ (see Figure~\ref{fig:UB}).
The scaling of the  energy is the same for the two variants of the construction.
Specifically, for any $k$, repeating the arguments above we introduce
\[
u_2^{L}(x_1,x_3):= ((-1)^k 2-\alpha x_1) x_3 + 2 \alpha x_1 k \delta +(-1)^{k+1} 2 k \delta,
\]
\[
u_2^{R}(x_1,x_3):= ((-1)^{k+1} 2-\alpha x_1) x_3 + 2 \alpha x_1 (k+1) \delta + (-1)^{k}  2 (k+1) \delta.
\]
The interface in the case of general $k$ is
\begin{equation}\label{interface2}
f_k^\alpha(x_1):=\frac{(-1)^k\alpha x_1 +1}{2}\delta  + k\delta.
\end{equation}
and \eqref{def_of_u2} now defines $u_2$ for both even and odd $k$. This way \eqref{el_bnd} holds in every $\Omega_k$, and $e_{23}$ is continuous across $x_3 = k$ interfaces.

To estimate the surface energy, we observe that $e_{23}(u_{\eps,\alpha})\in\{\pm 1\}$ almost everywhere implies that the jump is $\pm2$; the length of the jump set can be estimated by the total length of the oblique sides of the triangles. Each side has length $\sqrt{4+\delta^2}\le 2+\delta$, and since there are $2n$ of them
{and the thickness of the domain in the $x_2$ direction is $2$}
we obtain
\[
\Esur[u_{\eps,\alpha}] \le 2(4n+2n\delta)=8n+{4}.
\]
Choosing $n:=\lceil(\alpha^2/\eps)^{1/3}\rceil$ we obtain
\[
E_\ve[u_{\eps,\alpha}] \le c |\alpha|^{2/3}\eps^{2/3}+c\eps.
\]
If $\eps\le \alpha^2$ then $\eps\le \eps^{2/3}(\alpha^2)^{1/3}
=|\alpha|^{2/3}\eps^{2/3}$, the bound above
implies that $E_\ve[u_{\eps,\alpha}]\le c
|\alpha|^{2/3}\eps^{2/3}$ and concludes the proof {in the first case.}

{We now turn to the case $\alpha^2<\eps$,
{which is only relevant if $|\alpha|<1$}.}
We use a different construction, without microstructure. Precisely, we set
\begin{equation}\label{equbdefw}
 w(x):=\begin{pmatrix} 0 \\ 2x_3-\alpha x_1x_3 \\ \alpha x_1x_2\end{pmatrix}.
\end{equation}
Then
\begin{equation}
Dw(x)=\begin{pmatrix}
          0 & 0 & 0\\
          - \alpha x_3 & 0 & 2-\alpha x_1\\
          \alpha x_2 & \alpha x_1 &0
         \end{pmatrix}
\text{ and }
e(w)(x)=\begin{pmatrix}
          0 & -\frac12 \alpha x_3 & \frac12 \alpha x_2\\
          -\frac12 \alpha x_3 & 0 & 1\\
          \frac12 \alpha x_2 & 1&0
         \end{pmatrix}.
\end{equation}
Therefore $W(e(w))(x)\le \alpha^2$, and $\Esur[w]=0$. Therefore
\begin{equation}
 E_\eps[w]\le 8\alpha^2.
\end{equation}
If $\alpha^2\le \eps$ this concludes the proof.}
\end{proof}

We next address the upper bound for Theorem~\ref{theolowerboundneumann}.

{\begin{prop}\label{propubneu}
There is $c>0$ such that for all $\eps\in (0,\frac12]$ and all  $\gamma\in\R$ one can construct $u_{\eps,\gamma}\in  W^{1,2}(\Omega;\R^3)$ such that
\begin{equation*}
 E_\eps[u_{\eps,\gamma}]-\gamma M'(u_{\eps,\gamma})
  \le \min\{-\frac1c \gamma^2, c\eps^{2/3}-\frac1c|\gamma|\}.
\end{equation*}
\end{prop}}
\begin{proof}
 As discussed in {Section \ref{sec:getting-started}}, the two upper bounds arise from two very different constructions.

First we construct an affine deformation that generates the elastic response. This is similar to the function constructed in \eqref{equbdefw} for the regime of small $\alpha$, and is relevant in two regimes, either for very small forces (in the sense that $|\gamma|\lesssim\eps^{2/3}$)
or for very large ones (in the sense that $|\gamma|\gtrsim 1$).
Specifically, we set, for some $t\in\R$ chosen below,
\begin{equation}
 u_t(x):=\begin{pmatrix} tx_2x_3\\  2x_3 \\ 0\end{pmatrix}.
\end{equation}
Obviously $e_{23}(u_t)=1$ everywhere, for any $t$. Further, $E_\eps[u_0]=0$, and one can easily compute $E_\eps[u_t]=c_1 t^2$,
and $M'(u_t)=c_2 t$. Choosing $t=\frac{c_2}{2c_1}\gamma$, we obtain $E_\eps[u_t]-\gamma M'(u_t)=-\frac{c_2^2}{4c_1}\gamma^2$. This proves the first bound.

To prove the second one we use the deformation constructed in \eqref{eqdefuepsalpha} in the proof of Proposition~\ref{propubuniversal} for $\alpha=1$
if $\gamma<0$, and for $\alpha=-1$ if $\gamma\ge0$.
For simplicity of notation we focus on the second case.
As discussed above, one has
$E_\eps[u_{\eps,-1}]\le c\eps^{2/3}$.
Further, from \eqref{eqdefuepsalpha}
we obtain $u_1(x_1,x_2,1)-u_1(x_1,x_2,-1)=2x_2$, which implies $M'(u_{\eps,-1})=\frac83$. We obtain
$E_\eps[u_{\eps,-1}]-\gamma M'(u_{\eps,-1})\le c\eps^{2/3}-\frac83\gamma$, which concludes the proof.
\end{proof}

\section{Lower bounds} \label{sec:lower-bounds}
{Section \ref{subsec:LB-1} proves Proposition \ref{P1:rigid}, which asserts that $\alpha u^*$ is the
unique minimizer of the relaxed energy under any of our Dirichlet-type boundary conditions. As already
noted in Section \ref{sec:getting-started}, this result provides a conceptual foundation for all our
lower bound arguments. We then proceed, in Sections \ref{subsec:LB-2} -- \ref{subsec:LB-5}, to
prove the lower bounds for our Dirichlet-type boundary conditions. As already mentioned in
Section \ref{sec:getting-started}, we offer two different arguments. The first, presented
in Section \ref{subsec:LB-2}, uses a duality argument, taking advantage of the convexity
of the relaxed problem; it is restricted to $u \in \mathcal{A}_{TB} \cap \Aconstr$. The second
argument, presented in Sections \ref{subsec:LB-3}--\ref{subsec:LB-5}, is based on
a quantitative analogue of Proposition \ref{P1:rigid}; it works for all our Dirichlet-type boundary
conditions. Finally, Section \ref{subsec:LB-6} considers the functional associated with our
Neumann boundary condition. The lower bound for that functional relies again on the
quantitative analogue of Proposition \ref{P1:rigid}.
}

\subsection{Rigidity of exact solutions of the relaxed problem} \label{subsec:LB-1}
\begin{proof}[Proof of Proposition \ref{P1:rigid}]
One easily verifies that the function $u^*$ defined in (\ref{relaxed}) safisfies ${\Eel^\rel}[\alpha u^*] = 0$ {for all $\alpha\in[-1,1]$}. Further, it obeys the boundary conditions in
(\ref{BC alpha<1}) and, if $\alpha=1$, also those in
(\ref{BC alpha=1}). Therefore
$\alpha u^*\in \Athree\cap\Aconstr$ for all $\alpha\in[-1,1]$, and $u^*\in\Aone\cap\Aconstr$.

Assume now that $u$ is a function with
${\Eel^\rel}[u] = 0$. Then necessarily $e_{11} \equiv e_{22} \equiv e_{33} \equiv 0$ almost everywhere, which, in turn, implies that
\[
u(x)=\begin{pmatrix}
u_1(x_2,x_3)\\
u_2(x_1,x_3)\\
 u_3(x_1,x_2)
 \end{pmatrix}.
\]
Next, by $W(e(u))=0$ almost everywhere we obtain
\[
e_{13} \equiv 0 \Rightarrow \d_{3} u_1(x_2,x_3) \equiv - \d_{1}u_3(x_1,x_2).
\]
{The first term does not depend on $x_1$, the second does not depend on $x_3$. Therefore both depend only on $x_2$, and there is a function $\phi_2:(-1,1)\to\R$ such that}
\begin{equation}\label{eqrig13}
\partial_3u_1(x_2,x_3)=- \d_{1} u_3(x_1,x_2) = \phi_2(x_2),
\end{equation}
Similarly, swapping the indices 2 and 3,
\begin{equation}\label{eqrig12}
\partial_2u_1(x_2,x_3)=- \d_{1} u_2(x_1,x_3) = \phi_3(x_3).
\end{equation}
Taking the mixed derivative of these two conditions, we see that
\begin{equation}\label{eqrifds33}
 0=\partial_2\partial_3u_1-\partial_3\partial_2u_1=\partial_2\phi_2-\partial_3\phi_3
\end{equation}
distributionally, which implies
$\partial_2\phi_2(x_2)=\partial_3\phi_3(x_3)$ distributionally. Therefore both distributional derivatives are constant, and both $\phi_2$ and $\phi_3$ are affine.

Assume that $u\in \Aone$. Then $u_3(1,x_2)-u_3(-1,x_2)=2x_2$ and \eqref{eqrig13} give $\phi_2(x_2)=-x_2$ and $u_3=u^*_3$. Similarly,
\eqref{eqrig12} leads to
$\phi_3(x_3)=-x_3$ and $u_2=u^*_2$. Finally, $\partial_3 u_1(x_2,x_3)={-x_2}$ and $\partial_2 u_1(x_2,x_3)={-x_3}$ give
$u_1(x_2,x_3)={-x_2x_3}+c$, which concludes the proof
of \ref{P1:rigidone}.

Assume now that $u\in \Athree$. As $u_3$ does not depend on $x_3$, the boundary data immediately give $u_3=\alpha u_3^*$ everywhere, with $\phi_2(x_2)=-\alpha x_2$. Then \eqref{eqrifds33} gives $\phi_3(x_3)=-\alpha x_3+d$.
Integrating \eqref{eqrig13} and \eqref{eqrig12} leads to
\begin{equation}
 u_1(x_2,x_3)=-\alpha x_2x_3+dx_2+c
\end{equation}
and
\begin{equation}
 u_2(x_1,x_3)=\alpha x_1x_3-dx_1+\psi(x_3).
\end{equation}
Therefore any minimizer has the form given in \eqref{eqriglr} for some $c,d\in\R$ and some measurable function $\psi$.

Assume now that $u$ has the form given in
\eqref{eqriglr}.
One easily computes
$e_{ii}(u)=0$,
$e_{12}(u)=e_{13}(u)=0$ and
\begin{equation}
 2e_{23}(u)=2\alpha x_1+\psi'(x_3).
\end{equation}
Therefore $|e_{23}(u)|\le 1$ almost everywhere is equivalent to the fact that $\psi$ is Lipschitz with
$|\psi'|\le 2(1-|\alpha|)$ almost everywhere. This concludes the proof of \ref{P1:rigidthree}.
\end{proof}

\subsection{{Lower bound for $u \in \mathcal{A}_{TB} \cap \Aconstr$ using duality}} \label{subsec:LB-2}

{In this section, we use a duality argument to prove the lower bound when $u$ satisfies our top/bottom
boundary condition \eqref{BC alpha=1} and has the form \eqref{ansatz}.}

\begin{prop}\label{proplbdual}
There is $C>0$ such that for any
$\eps\in(0,1]$ and $u\in \Aconstr\cap\Aone$ one has
 \begin{equation}
  E_\eps[u]\ge C \eps^{2/3}.
 \end{equation}
\end{prop}

The key idea is that since the relaxed energy is convex, its convex dual can be used to bound it from below. The argument we
present was found by identifying the dual problem then using it with an appropriate class of test functions. But once the test functions
have been chosen, a duality-based lower bound proceeds by elementary calculations that use little more than integration by parts.
For maximum efficiency, we shall not discuss the dual problem; rather, we simply present the elementary calculations to which it led us.

The following duality-based lower bound does not use the constraint \eqref{ansatz}.

\begin{lem}\label{lemmax1aeel}
Let $u\in\Aone$.  For every $\eta\in W^{1,2}_0((-1,1)^2)$ and almost every $x_1^*\in(-1, 1)$ one has
\begin{equation}
\Eel[u]\ge \int_{(-1,1)^2}\left[
   (\partial_2 u_3-\partial_3 u_2)(x_1^*,\cdot,\cdot) \eta -\frac12(x_1^*+1) |{D}\eta|^2 \right]d\calL^2.
\end{equation}
\end{lem}
We remark that this bound, based on the boundary data and convexity, also holds for $\Eel^\rel$.
\begin{proof}
We shall prove the assertion for every $x_1^*$ such that the trace
$u(x_1^*,\cdot,\cdot)$ is in $W^{1,2}((-1,1)^2;\R^3)$.
By density it suffices to prove the assertion for $\eta\in C^\infty_c((-1,1)^2)$.
We first observe that for every $a,b\in\R$ one has
\begin{equation}
  2a^2\ge 2a b -\frac12b^2 .
\end{equation}
Using this with $a=e_{12}(u)(x)$ and $b=\partial_3\eta(x_2,x_3)$ and integrating
over $\Omega_*:=(-1,x_1^*)\times(-1,1)^2$
we obtain
\begin{equation}
 \int_{\Omega_*} 2e_{12}^2(u) dx \ge \int_{\Omega_*} 2e_{12}(u) \partial_3\eta d\calL^3 -\int_{\Omega_*} \frac12 (\partial_3\eta)^2 d\calL^3.
\end{equation}
Doing the same with $a'=e_{13}(u)(x)$ and $b'=-\partial_2\eta(x_2,x_3)$, summing
and using $W(e(u))\ge 2e_{12}^2(u)+2e_{13}^2(u)$
leads to
\begin{equation}\label{eqdualeepsnab}
\Eel[u]\ge \int_{\Omega_*} 2(e_{12}(u) \partial_3 \eta- e_{13}(u) \partial_2 \eta )d\calL^3 - \frac12(x_1^*+1)\int_{(-1,1)^2} |{D}\eta|^2
   d\calL^2.
\end{equation}
We now investigate the first integral in more detail. Writing it out explicitly, it is
\begin{equation}\begin{split}
            E^*:=&\int_{\Omega_*} 2(e_{12}(u) \partial_3 \eta- e_{13}(u) \partial_2 \eta )d\calL^3    \\
            =&\int_{\Omega_*} (\partial_2 u_1 \partial_3 \eta
            +    \partial_1 u_2 \partial_3 \eta
            -    \partial_1 u_3 \partial_2 \eta
            -    \partial_3 u_1 \partial_2 \eta)d\calL^3.
                \end{split}
\end{equation}
For almost every $x_1\in (-1,x_1^*)$ we have $u_1(x_1,\cdot,\cdot)\in W^{1,2}((-1,1)^2)$, and since $\eta\in C^\infty_c((-1,1)^2)$ we can integrate by parts the two terms involving $u_1$, leading to
%\begin{equation}
%\begin{split}
\begin{multline}
 \int_{(-1,1)^2} (\partial_2 u_1 \partial_3 \eta
            - \partial_3 u_1 \partial_2 \eta) d\calL^2(x_2,x_3)\\
            = - \int_{(-1,1)^2} u_1 (\partial_2\partial_3\eta-\partial_3\partial_2 \eta) d\calL^2(x_2,x_3)=0.
\end{multline}
%\end{split}
%\end{equation}
For the other two terms we use that $\eta$ does not depend on $x_1$, and obtain
\begin{equation}\label{eqlbdualestr2}
\begin{split}
            E^*
            =\int_{(-1,1)^2} ( u_2 \partial_3 \eta-
            u_3 \partial_2 \eta)
            (x_1^*,\cdot,\cdot)d\calL^2
            -\int_{(-1,1)^2}
            ( u_2 \partial_3 \eta-
            u_3 \partial_2 \eta)
            (-1,\cdot,\cdot)
            d\calL^2.
                \end{split}
\end{equation}
We next show that
the second integral vanishes, due to the boundary conditions \eqref{BC alpha=1}. Indeed, integrating by parts and recalling that $\eta\in C^\infty_c((-1,1)^2)$,
\begin{equation}\begin{split}
 \int_{(-1,1)^2}
            ( u_2 \partial_3 \eta-
            u_3 \partial_2 \eta)
            (-1,\cdot,\cdot)
            d\calL^2&=-
 \int_{(-1,1)^2}
           (x_3 \partial_3 \eta-
            x_2 \partial_2 \eta)
            (-1,\cdot,\cdot)
            d\calL^2\\
            &=
 \int_{(-1,1)^2}
           (\eta-
            \eta)
            d\calL^2=0.
            \end{split}
\end{equation}
In the first integral in \eqref{eqlbdualestr2} we can integrate by parts, since $u_{2,3}\in W^{1,2}((-1,1)^2)$, and obtain
\begin{equation}
 E^*            =-\int_{(-1,1)^2} ( \partial_3  u_2 -
            \partial_2  u_3 )
            (x_1^*,\cdot,\cdot)\eta d\calL^2.
\end{equation}
Recalling \eqref{eqdualeepsnab} and the definition of $E^*$ concludes the proof.
\end{proof}

The next Lemma shows that functions with small energy which obey the boundary conditions and the Ansatz \eqref{ansatz} obey a pointwise bound.
\begin{lem}\label{L:1N}
There exists a constant $C>0$
{such that for any $u\in\Aone\cap\Aconstr$, $\eps>0$, there is
a set $F\subseteq(-1,1)$ with $\calL^1((-1,1)\setminus F)\le \frac14$
such that
\begin{equation}\label{bound1}
|\hat u_3-x_1^*|(x_1^{*}, x_3) \leq C \Eel^{1/2}[u]
\end{equation}
for all
$x_1^*\in F$ and almost all $x_3 \in [-1,1]$.}
\end{lem}
\begin{proof}
We choose $x_2^*\in (\frac12,1)$ such that the trace of $u$ over $\{
x_2=x_2^*\}$ is in $W^{1,2}((-1,1)^2;\R^3)$ and
\begin{equation}\label{eqdefAl1}
 A:=\int_{(-1,1)^2} (e_{11}^2(u)+2e_{13}^2(u) +e_{33}^2(u))(x_1,x_2^*,x_3) dx_1dx_3 \le {2} \Eel[u].
\end{equation}
Let $w:=u-u_*$.
We observe that, from $e_{13}(u_*)=0$,
\begin{equation}
  2e_{13}(w)=\partial_1 w_3+\partial_3w_1
  = 2 e_{13}(u).
\end{equation}
Analogously $e_{ii}(w)=e_{ii}(u)$ for $i=1,3$, hence we can replace $w$ with $u$ in the definition of $A$ in \eqref{eqdefAl1}.
By Korn's inequality
applied to the map $w_{1,3}(\cdot, x_2^*,\cdot): (-1,1)^2\to\R^2$
there is $t\in\R$ such that
\begin{equation}\label{eqkornl1}
 \int_{(-1,1)^2} |\partial_1 w_3(\cdot, x_2^*,\cdot)-t|^2+
  |\partial_3 w_1(\cdot, x_2^*,\cdot)+t|^2 dx_1dx_3
  \le c A.
\end{equation}
From the boundary condition $w_3(\pm1,x_2,x_3)=0$ we obtain that
\begin{equation}
 \int_{(-1,1)^2} \partial_1 w_3 (\cdot, x_2^*,\cdot)dx_1dx_3=0,
\end{equation}
and with \eqref{eqkornl1} we obtain $|t|\le c A^{1/2}$. Therefore
\begin{equation}\label{eqkornl2}
 \int_{(-1,1)^2} |\partial_3 w_3(\cdot, x_2^*,\cdot)|^2+
 |\partial_1 w_3(\cdot, x_2^*,\cdot)|^2
 dx_1dx_3
 \le c A.
\end{equation}
We let $F$ be the set of those $x_1^*\in(-1,1)$ such that the trace of $w_3(\cdot, x_2^*,\cdot)$ obeys
\begin{equation}\label{eqp3w3}
 \int_{(-1,1)} |\partial_3 w_3(x_1^*, x_2^*,x_3)|^2dx_3
  \le 4c A,
\end{equation}
where $c$ is the same constant as in \eqref{eqkornl2}.
Obviously the measure of the complement is no larger than $1/4$. Fix any $x_1^*\in F$.
From \eqref{eqkornl2}, $w_3(\pm1,\cdot,\cdot)=0$ and the fundamental theorem of calculus applied in $x_1$ direction we also obtain
\begin{equation}\label{eqpx1w3}
 \int_{(-1,1)} |w_3(x_1^*,x_2^*,x_3)|^2dx_3 \le 2cA.
\end{equation}
Combining \eqref{eqp3w3} and \eqref{eqpx1w3} gives
\begin{equation}
  |w_3(x_1^*,x_2^*,x_3)|\le c A^{1/2} \text{ for almost every } x_3\in(-1,1).
\end{equation}
Recalling that $u\in \Aconstr$ implies $w_3(x)=x_2(\hat u_3(x_1,x_3)-x_1)$,
{\eqref{eqdefAl1},} and that $x_2^*\ge \frac12$ we conclude the proof.
\end{proof}

We finally come to the proof of the lower bound for functions in $\Aone\cap\Aconstr$.
\begin{proof}[Proof of Proposition~\ref{proplbdual}]
One key observation in this proof
is that
$u\in\Aconstr$ implies that $\partial_2u_3(x)=\hat u_3(x_1,x_3)$
depends only on $x_1$ and $x_3$, and analogously
$u_2$ (and hence $\partial_3 u_2$) also depends only on $x_1$ and $x_3$.  We define
$f:(-1,1)^2\to\R$ by
\begin{equation}
 f(x_1,x_3):=\hat u_3(x_1,x_3)+\partial_3u_2(x_1,x_3).
\end{equation}
Then $E_\eps[u]<\infty$ implies that $f\in BV((-1,1)^2;\{\pm 2\})$, with
$\Esur[u]=2\calH^1(J_f)$. In particular, for most choices of $x_1^*\in(-1,1)$ we have
$f(x_1^*,\cdot)\in BV((-1,1);\{\pm 2\})$ with  the number of jump points bounded by the energy. Precisely, we have
\begin{equation}\label{eqJf}
\# J_{f(x_1^*,\cdot)}\le 4 \Esur[u]
\end{equation}
outside an exceptional set of measure not larger than $1/8$.

We select $x_1^*\in (-\frac14,\frac14)$ such that \eqref{eqJf} as well as both
Lemma~\ref{lemmax1aeel} and Lemma~\ref{L:1N}
hold.
We define $\varphi\in L^2((-1,1))$ by
\begin{equation}\label{eqdefvarphi}
  \varphi(x_3):=\hat u_3(x_1^*,x_3)-\partial_3u_2(x_1^*,x_3).
\end{equation}
By Lemma~\ref{L:1N} we obtain
\begin{equation}
  |\hat u_3(x_1^*,x_3)-x_1^*|\le c \Eel^{1/2}[u]
  \text{ for almost all } x_3\in[-1,1].
\end{equation}
If $c \Eel^{1/2}[u]\ge\frac14$ then the proof is concluded. Therefore we are left with the case
\begin{equation}
  |\hat u_3(x_1^*,x_3)-x_1^*|\le \frac14  \text{ for almost all } x_3\in[-1,1].
\end{equation}
This implies
\begin{equation}\label{eqvarphif}
 |\varphi(x_3)+f(x_1^*,x_3)|
  =|2\hat u_3(x_1^*,x_3)|\le \frac12
+2|x_1^*| \le 1 \text{ for almost all } {x_3\in (-1,1)}.
\end{equation}
In particular, $|\varphi|\ge 1$ almost everywhere in $(-1,1)$, and it changes sign at most $N:=
\# J_{f(x_1^*,\cdot)}\le  4 \Esur[u]$ times.
For some $c_*\in\R$ we set
\begin{equation}\label{eqdefbeta}
 \beta(x_3):=\int_{-1}^{x_3} \varphi(t)dt -c_*
\end{equation}
and choose $c_*$ such that $\int_{(-1,1)} \beta(x_3)dx_3=0$.
Then $\beta'=\varphi$ yields
\begin{equation}\label{eqlbbeta2}
 \int_{(-1,1)} |\beta|^2(x_3)dx_3 \ge \frac{c}{(N+1)^2}\ge \frac{c}{(\Esur[u]+1)^2}.
\end{equation}
We define $\eta_0\in W^{1,2}_0((-1,1))$ by
\begin{equation}
 \eta_0(x_3):=\int_{-1}^{x_3} \beta(t)dt.
\end{equation}
For some $\theta\in C^1_c((-1,1))$ and $\gamma\in\R$ we set
\begin{equation}
 \eta(x_2,x_3):=
 \gamma\theta(x_2)\eta_0(x_3)=
 \gamma\theta(x_2)\int_{-1}^{x_3}\beta(t)dt.
\end{equation}
We observe that
$\eta\in W_0^{1,2}((-1,1)^2)$, with
$\partial_3 \eta(x_2,x_3)=\gamma\theta(x_2)\beta(x_3)$.
Lemma~\ref{lemmax1aeel} and \eqref{eqdefvarphi} then give
\begin{equation}\label{eqlbelb}
\Eel[u]\ge \int_{(-1,1)^2} \varphi(x_3)\eta(x_2,x_3) d\calL^2-{ \int_{(-1,1)^2}} |{D}\eta|^2 d\calL^2.
\end{equation}
We estimate
\begin{equation}
\int_{(-1,1)^2} |{D}\eta|^2 d\calL^2
\le \gamma^2 \|\theta\|_2^2 \|\beta\|_2^2 +
\gamma^2 \|\theta'\|_2^2\|\eta_0\|_2^2
\le C_\theta \gamma^2  \|\beta\|_2^2,
\end{equation}
where we used
$\|\eta_0\|_\infty\le \|\beta\|_1\le 2 \|\beta\|_2$, all norms being taken over $(-1,1)$. Further,
\begin{equation}
A:=\int_{(-1,1)^2} \varphi(x_3)\eta(x_2,x_3) d\calL^2
=\gamma \int_{(-1,1)} \theta(x_2)  dx_2 \int_{(-1,1)} \varphi(x_3)\eta_0(x_3) dx_3.
\end{equation}
Integrating by parts the second integral, recalling $\eta_0'=\beta$, $\eta_0(\pm1)=0$, and $\beta'=\varphi$, leads to
\begin{equation}
A
=\gamma \|\theta\|_1 \|\beta\|_2^2.
\end{equation}
We finally choose $\theta$ such that $\|\theta\|_1=1$ and then $\gamma:={1/2C_\theta}$, and obtain
from \eqref{eqlbelb}
\begin{equation}
 \Eel[u]\ge
 \gamma \|\beta\|_2^2
 - C_\theta\gamma^2\|\beta\|_2^2=
\frac{1}{{4C_\theta}} \|\beta\|_2^2.
\end{equation}
It only remains to combine this with \eqref{eqlbbeta2}. Indeed,
If $N=0$, then
\eqref{eqlbbeta2} gives $\|\beta\|_2^2\ge c$, hence $\Eel[u]\ge c\ge c \eps^{2/3}$. Otherwise,
\begin{equation}
 E_\eps[u]=
 \Eel[u]+\eps\Esur[u]
 \ge c \|\beta\|_2^2+\eps\frac{c}{\|\beta\|_2}
 \ge c\min_{t>0} (t^2+\frac\eps t)=c'\eps^{2/3}
\end{equation}
which concludes the proof.
\end{proof}

\subsection{{Lower bound for $u \in \mathcal{A}_{TB}$ or $u \in \mathcal{A}_{LR}^\alpha$ using a rigidity result}}
\label{subsec:LB-3}

{The following result is our rigidity-based lower bound for maps that obey our Dirichlet-type boundary conditions.
Its proof uses results that will be proved in the next two sections.}

\begin{prop}\label{proplowerbound}
 Let $\Omega:=(-1,1)^3$, $u\in W^{1,2}(\Omega;\R^3)$, $\eps>0$.
Then the following holds:
\begin{enumerate}
 \item\label{theolowerboundbcx1}
If $u\in\Aone$, which means
 \begin{equation}\label{eqbc}
 \begin{split}
  u_2(\pm1,x_2,x_3)=&\pm x_3,\\
  u_3(\pm1,x_2,x_3)=&\pm x_2,
 \end{split}
 \end{equation}
 then
\begin{equation}
 E_\eps[u]\ge c\min\{\eps^{2/3},1\}.
\end{equation}

 \item\label{theolowerboundbcx4}
If $u\in\Athree$ for some $\alpha\in[-1,1]$, which means
 \begin{equation}\label{eqbc3}
 \begin{split}
  u_3(x_1,x_2,\pm1)=& \alpha  x_1 x_2,
 \end{split}
 \end{equation}
 then
\begin{equation}
 E_\eps[u]\ge c\min\{|\alpha|^{2/3}\eps^{2/3},\alpha^2\}.
\end{equation}
\end{enumerate}
 \end{prop}

\begin{proof}[Proof of Proposition~\ref{proplowerbound}]
By Proposition~\ref{proplowerstructure}, we have
(for some $\bar s,\bar t,\beta\in\R$ and $b_2, d_3:(-1,1)\to\R$)
\begin{equation}\label{eqestu2bbcc}
\begin{split}
  \int_{\Omega} &
  \left[|u_2(x)-b_2(x_3)-\beta x_3 x_1 -\bar s x_1|^2
  \right.\\
  &+\left.
 |u_3(x)-d_3(x_2)-\beta x_1x_2 -\bar t x_1|^2
  \right]d\calL^3\le c E_\eps[u].
\end{split}
\end{equation}
In case \ref{theolowerboundbcx1} we have $\beta=1$,
{while}
in case
\ref{theolowerboundbcx4} we have $\beta=\alpha$.
By Proposition~\ref{proplocal} we then obtain
\begin{equation*}
E_\eps[u]\ge c \min\{\beta^2, \eps^{2/3}|\beta|^{2/3}\}
\end{equation*}
which concludes the proof.
\end{proof}

\subsection{Rigidity of approximate solutions of the relaxed problem} \label{subsec:LB-4}
We prove that any function $u$ with small elastic energy has approximately a specific form.
This result can be seen as a quantitative version of
Proposition \ref{P1:rigid}, and indeed the argument is similar to the
{one in Section~\ref{subsec:LB-1}}.
As in the case of
Proposition \ref{P1:rigid}, the rigidity estimate does not involve the surface energy, therefore we write the estimates in terms of $\Eel[u]\le E_\eps[u]$. The same estimates hold for the relaxed energy.

\begin{prop}\label{proplowerstructure}
Let $\Omega:=(-1,1)^3$, $u\in W^{1,2}(\Omega;\R^3)$.
Then there are $\beta\in\R$, two measurable functions $b_2,d_3:(-1,1)\to\R$, and  $\bar s,\bar t\in\R$, such that the estimates
\begin{equation}\label{eqestu2}
  \int_{\Omega}
  |u_2(x)-b_2(x_3)-\beta x_3 x_1 -\bar s x_1|^2
  d\calL^3 \le c \Eel[u],
\end{equation}
\begin{equation}\label{eqestu3}
  \int_{\Omega}
 |u_3(x)-d_3(x_2)-\beta x_1x_2 -\bar t x_1|^2
  d\calL^3 \le c \Eel[u],
\end{equation}
and
\begin{equation}\label{eqbdryM1}
 \int_{(-1,1)^2} \left| u_1(x_1,x_2,1)-u_1(x_1,x_2,-1)
 + 2\beta x_2+2\bar t\right| d\calL^2 \le c \Eel[u]^{1/2}
\end{equation}
hold. Further,
\begin{equation}\label{eqboundalpha}
 |\beta|\le 1 + c\Eel[u]^{1/3}+ c\Eel[u]^{1/2}.
\end{equation}
If the map $u$ obeys
\eqref{eqbc} then we can take $\beta=1$ in the previous estimates; if it obeys
\eqref{eqbc3}  then we can take $\beta=\alpha$.
 \end{prop}
The two conditions \eqref{eqestu2} and \eqref{eqestu3} characterize the behavior of $u_2$ and $u_3$ in the entire domain, and will be one key ingredient in the proof of the lower bound.  The boundary estimate \eqref{eqbdryM1} will instead be used to relate $\beta$ to the value of $M'(u)$
in Theorem~\ref{theolowerboundneumann}. We remark that other boundary estimates can also be obtained similarly, for other faces and other components; {but} we only state and prove explicitly the one that is used below in the proof of Theorem~\ref{theolowerboundneumann}.

We start by showing that if a function of two variables is close to being affine in one of them, and also close to being affine in the second one, then it is appoximately bilinear.
\begin{lem}\label{lemmau1rigid}
There is $c>0$ such that for any
functions $f,g,h,k\in L^2((-1,1))$, setting
\begin{equation*}
  \delta:=\int_{(-1,1)^2} |xf(y)+h(y)-yg(x)-k(x)|^2d\calL^2,
\end{equation*}
one can choose $\beta,F,G\in\R$ such that
\begin{equation*}
  \int_{(-1,1)} |f(y)-\beta y-F|^2 dy+
  \int_{(-1,1)}|g(x)-\beta x-G|^2 dx\le c\delta.
\end{equation*}
\end{lem}
\begin{proof}
We let $G$ be the average of $g$ and $K$ the average of $k$.
By convexity,
\begin{equation*}
  \int_{(-1,1)} |h(y)-yG-K|^2dy\le \delta,
\end{equation*}
so that with a triangular inequality
\begin{equation}\label{eqavfirst}
  \int_{(-1,1)^2} |xf(y)+y(G-g(x))+K-k(x)|^2d\calL^2\le 4\delta.
\end{equation}
Averaging this expression in $y$, and letting $F$ be the average of $f$,
\begin{equation*}
  \int_{(-1,1)} |xF+K-k(x)|^2dx\le 4\delta
\end{equation*}
so that with another triangular inequality \eqref{eqavfirst} gives
\begin{equation}\label{eqsyfFgG}
  \int_{(-1,1)^2} |x(f(y)-F)-y(g(x)-G)|^2d\calL^2\le c\delta.
\end{equation}
We define $\beta$ as  the average of $y\mapsto (f(y)-F)/y$ over $(-1,1)\setminus (-\frac12,\frac12)$,
\begin{equation}
 \beta:=\int_{(-1,-\frac12)\cup(\frac12,1)}
 \frac{f(y)-F}{y} dy.
\end{equation}
Then a similar convexity argument leads to
\begin{equation*}
\begin{split}
&\int_{(-1,1)} |x\beta-(g(x)-G)|^2dx\\
 & \le
\int_{(-1,1)^2\setminus (-1,1)\times (-\frac12,\frac12)} \left|x\frac{f(y)-F}{y}-(g(x)-G)\right|^2d\calL^2\\
& \le
4\int_{(-1,1)^2} \left|x(f(y)-F)-y(g(x)-G)\right|^2d\calL^2\le c\delta,
\end{split}
\end{equation*}
where in the last step we used \eqref{eqsyfFgG}.
Analogously, letting $\eta$ be the average of $(g(x)-G)/x$ over the same set, and using again \eqref{eqsyfFgG},
\begin{equation*}
  \int_{(-1,1)} |y\eta-(f(y)-F)|^2dy
\le c\delta,
\end{equation*}
Inserting in \eqref{eqsyfFgG} we obtain $|\beta-\eta|^2\le c\delta$. A triangular inequality concludes the proof.
\end{proof}

\begin{proof}[Proof of Proposition~\ref{proplowerstructure}]
To simplify notation we write $E:=\Eel[u]$; we can assume that $E<\infty$, which is the same as $|e_{23}(u)|=1$ almost everywhere.

\emph{Step 1: Korn's inequality on slices.}
For a fixed $x_3\in(-1,1)$, we consider the slice $v^{(x_3)}:(-1,1)^2\to\R^2$, $v^{(x_3)}(x_1,x_2):=(u_1,u_2)(x_1,x_2,x_3)$.
For almost every $x_3$, we have $v^{(x_3)}\in W^{1,2}((-1,1)^2;\R^2)$, and ${D}  v^{(x_3)}(x_1,x_2)$ coincides almost everywhere with the matrix obtained
dropping the third row and the third column of
${D} u(x)$.
By Korn's inequality and Poincar\'e's inequality there is an affine isometry $A^{(x_3)}:\R^2\to\R^2$, of the form
\begin{equation*}
  A^{(x_3)}(x_1,x_2)=\begin{pmatrix} b_1(x_3)-s(x_3)x_2\\b_2(x_3)+s(x_3)x_1 \end{pmatrix}
\end{equation*}
for some measurable $b:(-1,1)\to\R^2$ and $s:(-1,1)\to\R$,
such that
\begin{equation}\label{eqvx3w12}
\begin{split}
 \int_{(-1,1)^2} |v^{(x_3)}-A^{(x_3)}|^2
 &+|{D} v^{(x_3)}-{D} A^{(x_3)}|^2
 d\calL^2 \\
 &\le c \int_{(-1,1)^2} |e'(v^{(x_3)})|^2 d\calL^2 = c e^{(3)}(x_3),\\
\end{split}
\end{equation}
where
\begin{equation*}
e^{(3)}(x_3):=\int_{(-1,1)^2} (e_{11}^2+2e_{12}^2+e_{22}^2)(u)(\cdot,\cdot, x_3)d\calL^2.
\end{equation*}
By the trace theorem,
\begin{equation}\label{eqvxtrace}
\begin{split}
 &\int_{\partial(-1,1)^2} |v^{(x_3)}-A^{(x_3)}|^2
 d\calH^1\le ce^{(3)}(x_3).
\end{split}
\end{equation}
We integrate \eqref{eqvx3w12} over $x_3$, dropping the second term and inserting
the explicit form of $v^{(x_3)}$ and $A^{(x_3)}$, and obtain
\begin{equation}\label{equ13omega}
  \int_{\Omega} |u_1(x)-b_1(x_3)+s(x_3)x_2|^2
  +|u_2(x)-b_2(x_3)-s(x_3)x_1|^2
  d\calL^3 \le c E
\end{equation}
and, proceeding similarly from \eqref{eqvxtrace},
\begin{equation}\label{equ13pomega}
  \int_{\partial_3\Omega} |u_1(x)-b_1(x_3)+s(x_3)x_2|^2
  +|u_2(x)-b_2(x_3)-s(x_3)x_1|^2
  d\calH^2 \le c E.
\end{equation}
Here $\partial_3\Omega:=\{x: (x_1,x_2)\in\partial(-1,1)^2,x_3\in(-1,1)\}=\partial\Omega\cap\{|x_3|<1\}$ is the part of the boundary of $\Omega$ whose normal is not $\pm e_3$.

The same argument can be performed swapping $x_2$ and $x_3$.
For $x_2\in(-1,1)$, we define $w^{(x_2)}\in W^{1,2}((-1,1)^2;\R^2)$  by $w^{(x_2)}(x_1,x_3):=(u_1,u_3)(x_1,x_2,x_3)$.
By the Korn-Poincar\'e inequality there is an affine isometry $B^{(x_2)}:\R^2\to\R^2$, of the form
\begin{equation*}
B^{(x_2)}(x_1,x_3)=\begin{pmatrix} d_1(x_2)-t(x_2)x_3\\d_3(x_2)+t(x_2)x_1 \end{pmatrix}
\end{equation*}
for some measurable $d:(-1,1)\to\R^2$ and $t:(-1,1)\to\R$,
such that
\begin{equation*}%\label{eqvx3w13}
\begin{split}
 \int_{(-1,1)^2} |w^{(x_2)}-B^{(x_2)}|^2
 &+|{D} w^{(x_2)}-{D} B^{(x_2)}|^2
 d\calL^2
\le c \int_{(-1,1)^2} |e(w^{(x_2)})|^2 d\calL^2 .
\end{split}
\end{equation*}
The same argument as above leads to
\begin{equation}\label{equ12omega}
  \int_{\Omega} |u_1(x)-d_1(x_2)+t(x_2)x_3|^2
  +|u_3(x)-d_3(x_2)-t(x_2)x_1|^2
  d\calL^3 \le c E
\end{equation}
and
\begin{equation}\label{equ12pomega}
  \int_{\partial_2\Omega} |u_1(x)-d_1(x_2)+t(x_2)x_3|^2
  +|u_3(x)-d_3(x_2)-t(x_2)x_1|^2
  d\calH^2 \le c E,
\end{equation}
with $\partial_2\Omega:=\{x: (x_1,x_3)\in\partial(-1,1)^2,x_2\in(-1,1)\}=\partial\Omega\cap\{|x_2|<1\}$.

\emph{Step 2: Structure of the functions $t$ and $s$.}
The two volume estimates permit, via Lemma~\ref{lemmau1rigid}, to prove that $t$ and $s$ are approximately affine.
The key observation is that the component $u_1$ is estimated both  in the first term of \eqref{equ13omega} and in the first one of \eqref{equ12omega}, therefore a triangular inequality gives
\begin{equation*}
  \int_{\Omega} |b_1(x_3)-s(x_3)x_2-d_1(x_2)+t(x_2)x_3|^2
  d\calL^3 \le c E.
\end{equation*}
The integrand does not depend on $x_1$, hence the integral is effectively only over $(-1,1)^2$.
Using Lemma~\ref{lemmau1rigid}
shows that there are $\beta$, $\bar s$, $\bar t\in\R$ such that
\begin{equation}\label{eqalpharigid}
  \int_{(-1,1)} |s(x_3)-\beta x_3-\bar s|^2  dx_3+
  \int_{(-1,1)} |t(x_2)-\beta x_2-\bar t|^2  dx_2
 \le c E.
 \end{equation}
Using this estimate in
the second term of \eqref{equ13omega} and the second one of \eqref{equ12omega}, with a triangular inequality one immediately obtains
\eqref{eqestu2} and \eqref{eqestu3}.
Inserting in \eqref{equ12pomega} leads to
\begin{equation}\label{equ12pomegab}
  \int_{\partial_2\Omega} |u_1(x)-d_1(x_2)+\beta x_2x_3+\bar t x_3|^2
  +|u_3(x)-d_3(x_2)-\beta x_1x_2-\bar t x_1|^2
  d\calH^2 \le c E.
\end{equation}

\emph{Step 3: Boundary terms.}
We start from the boundary estimate in
\eqref{equ12pomegab}. As $(-1,1)^2\times\{-1,1\}\subseteq\partial_2\Omega$,
from the first term we obtain
\begin{equation*}\begin{split}
 & \int_{(-1,1)^2} \Bigl[|u_1(x_1,x_2,1)-d_1(x_2)+\beta x_2+\bar t|^2 \\
& \hskip1cm + |u_1(x_1,x_2,-1)-d_1(x_2)-\beta x_2-\bar t|^2 \Bigr]
  d\calL^2 \le c E,
\end{split}
\end{equation*}
and
with a triangular inequality
\begin{equation*}
 \int_{(-1,1)^2} \left| u_1(x_1,x_2,1)-u_1(x_1,x_2,-1)
 + 2\beta x_2+2\bar t\right|^2d\calL^2  \le c E
\end{equation*}
which implies \eqref{eqbdryM1}.

Assume now that
\eqref{eqbc3} holds.
Using again $(-1,1)^2\times\{-1,1\}\subseteq\partial_2\Omega$,
from the second term of \eqref{equ12pomegab} we obtain
\begin{equation*}
  \int_{(-1,1)^2}
  |\alpha x_1x_2-d_3(x_2)-\beta x_2x_1-\bar t x_1|^2
  d\calL^2 \le c E,
\end{equation*}
which leads to $|\alpha-\beta|\le c E^{1/2}$.

Finally, assume that \eqref{eqbc} holds.
Using $\{-1,1\}\times(-1,1)^2\subseteq\partial_2\Omega$,
from the second term of \eqref{equ12pomegab} we obtain
\begin{equation*}
  \int_{(-1,1)^2} |x_2-d_3(x_2)-\beta x_2-\bar t |^2
  + |-x_2-d_3(x_2)+\beta x_2+\bar t |^2
  d\calL^2 \le c E.
\end{equation*}
With a triangular inequality
we obtain
\begin{equation*}
  \int_{(-1,1)} |2(1-\beta)x_2-2\bar t |^2
    dx_2 \le c E
\end{equation*}
and therefore $|\beta-1|\le c E^{1/2}$.

\emph{Step 4: Upper bound on $\beta$.}
It remains to prove the estimate \eqref{eqboundalpha}.
For any $\varphi\in C^1_c(\Omega)$, from $|\partial_2 u_3+\partial_3u_2|=2$ almost everywhere we obtain
\begin{equation}\label{equpartphi}
 \left|\int_\Omega (u_3\partial_2 \varphi + u_2 \partial_3 \varphi) d\calL^3\right|
=\left| \int_\Omega \varphi (\partial_2 u_3+\partial_3 u_2) d\calL^3\right| \le 2\int_\Omega|\varphi|d\calL^3.
\end{equation}
We define
\begin{equation*}
 G:=\int_\Omega \left[(d_3(x_2)+\beta x_1x_2+\bar t x_1)\partial_2\varphi+
 (b_2(x_3)+\beta x_3x_1+\bar s x_1)\partial_3\varphi\right]
 d\calL^3 .
\end{equation*}
Using \eqref{eqestu2}, \eqref{eqestu3} and Hölder's inequality, with ${D}'\varphi:=(\partial_2\varphi,\partial_3\varphi)$,
\begin{equation*}
\left| \int_\Omega (u_3\partial_2 \varphi + u_2 \partial_3 \varphi )d\calL^3-G \right|
\le c \|{D}'\varphi\|_{L^2(\Omega)} E^{1/2},
\end{equation*}
and with \eqref{equpartphi} we obtain
\begin{equation}\label{eqGl2}
 |G|\le
 2\int_\Omega|\varphi|d\calL^3+
 c \|{D}'\varphi\|_{L^2(\Omega)} E^{1/2}.
\end{equation}
At the same time, using that
\begin{equation*}
\int_{(-1,1)} \partial_2\varphi(x_1,x_2',x_3) dx_2'=
\int_{(-1,1)} \partial_3\varphi(x_1,x_2,x_3') dx_3'=0,
\end{equation*}
\begin{equation*}\begin{split}
 G =&\int_\Omega \left[d_3(x_2)
 \partial_2\varphi + b_2(x_3)\partial_3\varphi
 +\beta x_1x_2 \partial_2\varphi
 +\beta x_1x_3 \partial_3\varphi \right]
  d\calL^3 .
\end{split}
\end{equation*}
We assume now that $\varphi(x_1,x_2,x_3)=-\varphi(-x_1,x_2,x_3)$ for all $x\in\Omega$,
which implies the same symmetry for $\partial_2\varphi$ and $\partial_3 \varphi$. Then the first two terms in the last integral disappear. We integrate by parts the remaining ones and conclude
\begin{equation}\label{eqGalpha}
\begin{split}
 G=&
-2\beta
 \int_\Omega  x_1\varphi
 d\calL^3.
\end{split}
\end{equation}
Combining \eqref{eqGalpha} and \eqref{eqGl2}, we conclude that
\begin{equation}\label{eqalphaGaphi}
2 |\beta|
\left| \int_\Omega  x_1\varphi
 d\calL^3\right| =
 |G|\le c \|{D}' \varphi\|_{L^2(\Omega)} E^{1/2}
 + 2\int_\Omega|\varphi|d\calL^3
\end{equation}
for every $\varphi\in C^1_c(\Omega)$ such that $\varphi(x_1,x_2,x_3)=-\varphi(-x_1,x_2,x_3)$.

We next choose the function $\varphi$.
We fix $\psi\in C^1_c((-1,1)^2;[0,\infty))$
with $\|\psi\|_{L^1((-1,1)^2)}=1$ and set $\varphi(x):=(\theta(x_1)-\theta(-x_1))\psi(x_2,x_3)$,
where
$\theta\in C^1_c((0,1);[0,\infty))$ will be chosen below.
Then
${D}'\varphi(x)=(\theta(x_1)-\theta(-x_1)){D}\psi(x_2,x_3)$ and
\begin{equation*}
\|{D}'\varphi\|_{L^2(\Omega)}=\sqrt2 \|\theta\|_{L^2((0,1))} \|{D} \psi\|_{L^2((-1,1)^2)}.
\end{equation*}
Inserting these expressions for $\varphi$ in \eqref{eqalphaGaphi} then leads to
\begin{equation*}
4 |\beta|
 \int_{(0,1)} x_1\theta(x_1) dx_1 \le C
 \|\theta\|_{L^2((0,1))}
  E^{1/2}
 + 4\int_{(0,1)} \theta(x_1) dx_1
\end{equation*}
for all $\theta\in C^1_c((0,1);[0,\infty))$, with a constant
$C$ that may depend on $\psi$, but not on $\theta$, $\eps$ or $u$. By density the same holds for
any $\theta\in L^2((0,1);[0,\infty))$. We select $\theta:=\chi_{(1-\delta,1)}$ for some $\delta\in (0,\frac12]$, and obtain
\begin{equation*}
 4|\beta| (1-\delta)\delta \le C \delta^{1/2}   E^{1/2}
+4\delta,
\end{equation*}
which implies
\begin{equation*}
 |\beta| \le C  \frac{  E^{1/2}}{\delta^{1/2}}
+\frac 1{1-\delta}.
\end{equation*}
If $E\le 1$ then
selecting
$\delta:=\frac12 E^{1/3}$ concludes the proof of
\eqref{eqboundalpha}.
If instead
 $E> 1$ we consider $\delta:=1/2$.
\end{proof}

\subsection{Auxiliary lower bound for approximately bilinear deformations.} \label{subsec:LB-5}
In this Section we prove that if a function is approximately bilinear, in the sense made precise in
\eqref{eqestu2bb} below, then the energy cannot be small. As the example $u^*$ shows, this cannot be obtained from the elastic energy alone. One key ingredient is a rigidity result, presented in Lemma~\ref{lemmaccf} below, which shows that
functions with finite energy are appoximately affine on a scale set by the surface energy. Therefore, either the energy is large, or the function is not so close to the quadratic expression, and we obtain a lower bound of the form \eqref{eqlbefeps}.
\begin{prop}\label{proplocal}
  Let $\Omega:=(-1,1)^3$, $u\in W^{1,2}(\Omega;\R^3)$, $\eps>0$.
For some $\beta\in\R$,  $b_2,d_3:(-1,1)\to\R$ measurable, $\bar s,\bar t\in\R$ let
\begin{equation}\label{eqestu2bb}
\begin{split}
  F:=\int_{\Omega} &
  \left[|u_2(x)-b_2(x_3)-\beta x_3 x_1 -\bar s x_1|^2
  \right.\\
  &+\left.
 |u_3(x)-d_3(x_2)-\beta x_1x_2 -\bar t x_1|^2
  \right]d\calL^3.
\end{split}
\end{equation}
Assume $|\beta|\le 3$. Then
\begin{equation}\label{eqlbefeps}
F+ E_\eps[u]\ge c \min\{\beta^2, \eps^{2/3}|\beta|^{2/3}\}.
\end{equation}
\end{prop}
Before starting the proof,
let us sketch the main strategy. The key idea is common to many lower bounds for singularly perturbed
nonconvex problems: on a suitable length scale, called $\lambda$ below, either the surface energy is
large, or one of the phases dominates. In the latter case, the function necessarily deviates from the
relaxed solution. In order to make this precise,
fix $\eta\ll1$ and
consider a typical cube $Q_\lambda\subset\Omega$ of side $\lambda$.
If the part of the boundary of $\{e_{23}(u)=1\}$ inside $Q_\lambda$ is larger than $\eta\lambda^2$, then
the total length of the boundary is at least $\eta\lambda^2/\lambda^3=\eta/\lambda$. If instead it is
smaller than $\eta\lambda^2$, then one of the two phases dominates. Assume it is $A:=Q_\lambda\cap \{ e_{23}(u)=1\}$.
In order to learn that $u$ is approximately affine in this cube, we need to use Korn's inequality on the set $A$;
in order to obtain the optimal scaling of the lower bound the  constant cannot depend on $A$. However, $A$ is a
set of finite perimeter, and might be very irregular. This difficulty is solved resorting to the Korn-Poincar\'e
inequality with holes presented in Lemma~\ref{lemmaccf} below. One final twist of the proof is that we need to
contrast the assumption that $u$ 
{is close to a bilinear function in the sense of \eqref{eqestu2bb}}, hence we need to take two cubes and consider the difference in
behavior among the two cubes, see \eqref{eqdefvfrow} and following arguments.

We start recalling the following special case of \cite[Th.~1.1]{ChambolleContiFrancfort2016}, which gives the
extension of the Korn-Poincar\'e inequality needed to obtain local rigidity.
\begin{lem}\label{lemmaccf}
Let $\Omega\subseteq\R^n$ be a bounded connected Lipschitz set. There is $c>0$ such that for any $u\in SBD^2(\Omega)$
there are an affine linear isometry $a:\R^n\to\R^n$ and a Borel set $\omega\subseteq \Omega$ such that
\begin{equation}\label{lemmaccf1}
 \|u-a\|_{L^2(\Omega\setminus \omega)}\le c  \|e(u)\|_{L^2(\Omega)}
\end{equation}
and
\begin{equation}\label{lemmaccf2}
 \mathcal L^n(\omega)\le c ({\mathcal H^{n-1}}(J_u) )^{n/(n-1)}.
\end{equation}
\end{lem}
Here and below, an affine map $a:\R^n\to\R^n$ is a linear isometry if $Da+Da^T=0$. We recall that $BD(\Omega)$ is
the set of $u\in L^1(\Omega;\R^n)$ such that the distributional strain $Eu:=\frac12(Du+Du^T)$ is a bounded measure;
one can prove that $Eu=e(u)\calL^n + E^cu+[u]\odot \nu \calH^{n-1} \LL J_u$, with $J_u$ the $n-1$-rectifiable jump
set of $u$, $\nu$ the normal and  $[u]$ the jump,
and $E^cu$ orthogonal to $\calL^n$ and vanishing on sets of finite $n-1$-dimensional measure. Further, $SBD^2(\Omega)$
is the set of those $u\in BD(\Omega)$ with $E^cu=0$, $e(u)\in L^2(\Omega;\R^{n\times n}_\sym)$ and $\calH^{n-1}(J_u)<\infty$.
In particular, if $u,v\in W^{1,2}(\Omega;\R^n)$ and $\borelset\subseteq\Omega$ is a set of finite perimeter, then the
function $w:=u\chi_\borelset + v\chi_{\Omega\setminus \borelset}$ is in $SBD^2(\Omega)$,
with $J_w\subset \Omega\cap \partial^* \borelset$ and
$\nabla w=(\nabla u)\chi_\borelset + (\nabla v)\chi_{\Omega\setminus \borelset}$. Here and below, $\partial^*$
denotes the measure-theoretic boundary of a set. We refer to \cite{AmbrosioCosciaDalmaso1997} for standard properties
of functions of bounded deformation.

We shall use the following corollary of Lemma \ref{lemmaccf}:
\begin{lem}\label{lemmaccfcube}
Let $n\ge 1$, $\gamma\in (0,1]$.
There are $\eta\in(0,1]$ and $c>0$ such that for any
cube
$Q_r\subset\R^n$ of side $r$ the following holds:
\begin{enumerate}
\item\label{lemmaccfcuberigid}
For any $u\in SBD^2(Q_r)$ with $\calH^{n-1}(J_u)\le\eta r^{n-1}$ there are an affine linear isometry $a:\R^n\to\R^n$ and a Borel set $\omega\subseteq Q_r$ such that
\begin{equation}
 \|u-a\|_{L^2(Q_r\setminus \omega)}\le c r \|e(u)\|_{L^2(Q_r)}
\end{equation}
and
\begin{equation}
 \mathcal L^n(\omega)\le \gamma r^n.
\end{equation}
\item\label{lemmaccfcubeisop}
For any Borel set $\borelset\subseteq Q_r$ with $\calH^{n-1}(Q_r\cap \partial^* \borelset)\le\eta r^{n-1}$ one has
\begin{equation}
 \min\{\calL^n(\borelset), \calL^n(Q_r\setminus \borelset)\} \le \gamma r^n.
\end{equation}
\end{enumerate}
\end{lem}
\begin{proof}
By scaling it suffices to prove the assertion for $Q_1=(0,1)^n$. Let $c_1$ be the constant in
Lemma \ref{lemmaccf} for $\Omega=Q_1$.
Then \eqref{lemmaccf2} implies
\begin{equation*}
 \mathcal L^n(\omega)\le c_1 ({\mathcal H^{n-1}}(J_u) )^{n/(n-1)}
 \le c_1\eta^{n/(n-1)}
\end{equation*}
so that the first assertion holds with $c:=c_1$ and any $\eta$ such that $c_1\eta^{n/(n-1)}\le\gamma$.

The second assertion follows from the standard relative isoperimetric inequality or, equivalently, from the
Poincar\'e inequality for the characteristic function of $\borelset$:
\begin{equation*}
\frac12 \min\{\calL^n(\borelset), \calL^n(Q_1\setminus \borelset)\} \le
\inf_{b\in\R}\|\chi_\borelset-b\|_{L^1(Q_1)} \le c_P |D\chi_\borelset|(Q_1) \le c_P\eta,
\end{equation*}
so that the assertion holds for any $\eta\le \gamma/(2c_P)$.
\end{proof}

The second ingredient in the proof of Proposition~\ref{proplocal} is a method to transform $L^p$ estimates on
second-degree polynomials on large subsets of a cube into estimates on the coefficients. To keep notation
simple we only discuss the specific version used below.
\begin{lem}\label{lemmacubelin}
There is $c>0$ such that for any measurable
$\omega\subseteq Q_r:=x^*+(-\frac12r,\frac12r)^3$ with $\calL^3(\omega)\le \frac14r^3$ and
any $z\in\R^3$ one has
\begin{equation*}
 r^2|z|^2
  \le c \min_{q\in\R} \frac1{r^3}\int_{Q_r\setminus\omega} |z\cdot x - q|^2 dx.
\end{equation*}
\end{lem}
\begin{proof}
We only prove the bound on $z_1$.
Let $T_1:\R^3\to\R^3$ be the reflection along $e_1$ which leaves $Q_r$ invariant,
$T_1(x):=(2x_1^*-x_1,x_2,x_3)$. Let $\tilde\omega:=\omega\cup T_1(\omega)$ and
$g(x):=z\cdot x-q$. From
\begin{equation*}
 \|g\circ T_1-g\|_{L^2(Q_r\setminus\tilde\omega)}
 \le
 \|g\circ T_1\|_{L^2(Q_r\setminus\tilde\omega)}+
 \|g\|_{L^2(Q_r\setminus\tilde\omega)}
\end{equation*}
we obtain, as $T_1(x)-x=(2(x_1^*-x_1),0,0)$,
\begin{equation*}
 2\|z_1(x_1^*-x_1)\|_{L^2(Q_r\setminus\tilde\omega)}\le 2 \|g\|_{L^2(Q_r\setminus\omega)}.
\end{equation*}
Consider now the set
\begin{equation*}
A:=(Q_r\setminus \tilde\omega)
\cap\{|x_1-x_1^*|\ge \frac1{8} r\}.
\end{equation*}
From $\calL^3(\tilde\omega)\le 2\calL^3(\omega)\le\frac12r^3$ we
obtain $\calL^3(A)\ge r^3(1-\frac12-\frac14)=\frac14r^3$.
Therefore
\begin{equation*}
 2 |z_1| \frac{r}{8} (\calL^3(A))^{1/2} \le 2  \|g\|_{L^2(Q_r\setminus\omega)}
\end{equation*}
which concludes the proof.
\end{proof}

We finally come to the proof of Proposition~\ref{proplocal}.
\begin{proof}[Proof of Proposition~\ref{proplocal}]
To simplify notation we set $E:=F+E_\eps[u]$. We can assume $E<\infty$.
Fix $\lambda\in(0,\frac14]$, chosen below.
We consider the measure
\begin{equation*}
\begin{split}
 \mu:= W(e(u))\calL^3+\eps |De_{23}(u)| + &
   |u_2(x)-b_2(x_3)-\beta x_3 x_1 -\bar s x_1|^2\calL^3 \\
   &+
   |u_3(x)-d_3(x_2)-\beta x_1x_2 -\bar t x_1|^2\calL^3.
\end{split}
\end{equation*}
By \eqref{eqestu2bb}, $\mu(\Omega)= E$.
First we pick $(x^*_2,x^*_3)\in (-1+\frac\lambda2,1-\frac\lambda2)^2$ such that
\begin{equation*}
 \mu( (-1,1)\times (x^*_2-\frac12\lambda,x^*_2+\frac12\lambda)
 \times (x^*_3-\frac12\lambda,x^*_3+\frac12\lambda)) \le \lambda^2 \mu(\Omega).
\end{equation*}
Then we pick $x_1^*\in ( -\frac12,0)$ such that, setting $y^*:=(x_1^*+\frac12,x_2^*,x_3^*)$,
the two disjoint cubes  $Q_\lambda:= x^*+(-\frac12\lambda,\frac12\lambda)^3\subseteq\Omega$ and
$\hat Q_\lambda:= y^*+(-\frac12\lambda,\frac12\lambda)^3=Q_\lambda+\frac12 e_1\subseteq\Omega$
have the property
\begin{equation*}
 \mu(Q_\lambda\cup \hat Q_\lambda) \le c \lambda^3 \mu(\Omega).
\end{equation*}
In particular, this implies
\begin{equation}\label{eqchoiceEb}
E_\eps[u,Q_\lambda \cup \hat Q_\lambda]
\le c\lambda^3E
\end{equation}
and
\begin{equation}\label{eqchoicel2b}
\begin{split}
    \int_{Q_\lambda\cup \hat Q_\lambda}
  |u_2(x)-b_2(x_3)-\beta x_3 x_1 -\bar s x_1|^2
  d\calL^3 \le c   \lambda^3  E,
 \\
 \int_{Q_\lambda\cup \hat Q_\lambda}
 |u_3(x)-d_3(x_2)-\beta x_1x_2 -\bar t x_1|^2
  d\calL^3 \le c   \lambda^3  E.
\end{split}
\end{equation}
Let $\eta>0$ be as in Lemma \ref{lemmaccfcube} with $\gamma=\frac1{12}$. Now distinguish two cases.
If $|De_{23}(u)|(Q_\lambda\cup\hat Q_\lambda)\ge \eta \lambda^2$ then \eqref{eqchoiceEb} gives
\begin{equation}\label{eqesurfenb}
E\ge c \lambda^{-3} E_\eps[u,Q_\lambda\cup\hat Q_\lambda] \ge c \lambda^{-3} \eps \eta \lambda^2 = c \eta \frac{\eps}{\lambda}.
\end{equation}
If instead $|De_{23}(u)|(Q_\lambda\cup\hat Q_\lambda)< \eta \lambda^2$, we
let $f:=e_{23}(u)\in SBV(\Omega;\{-1,1\})$ and
define
$w:\Omega\to\R^3$ by
\begin{equation*}
 w(x):=
 u(x)-2x_3 f(x) e_2=
 \begin{cases}
        u(x) - 2x_3e_2, & \text{ if } e_{23}(u)(x)=1,\\
        u(x) + 2x_3e_2, & \text{ otherwise.}\\
       \end{cases}
\end{equation*}
Since $u\in W^{1,2}(\Omega;\R^3)$ and $f\in SBV(\Omega;\{\pm1\})$, we have
$w\in SBV^2(\Omega;\R^3 )\subseteq SBD^2(\Omega)$, with $J_w\subseteq J_f$ up to $\calH^2$-null sets.
We next consider the part of the strain which is absolutely continuous with respect to the Lebesgue measure. Since $\nabla f=0$,
\begin{equation}
e_{23}(w)=e_{23}(u)-f=0
\end{equation}
and
\begin{equation*}
 |e|^2(w)=(e_{11}^2+e_{22}^2+e_{33}^2+2e_{12}^2+2e_{13}^2)(u)\le W(e(u))
\end{equation*}
almost everywhere, so that
\begin{equation*}
\int_{Q_\lambda\cup \hat Q_\lambda} |e|^2(w) d\calL^3\le E_\eps[u,Q_\lambda\cup \hat Q_\lambda].
\end{equation*}
We next apply Lemma \ref{lemmaccfcube}\ref{lemmaccfcubeisop}
to the set $Q_\lambda\cap\{f=1\}$,
and then the same on $\hat Q_\lambda$. This is admissible since $|Df|(Q_\lambda
\cup \hat Q_\lambda)<\eta\lambda^2$.
We obtain that there are $\sigma\in\{-1,1\}$ and $\hat\sigma\in\{-1,1\}$, such that
\begin{equation*}
\calL^3(Q_\lambda\cap\{f\ne \sigma\})+ \calL^3(\hat Q_\lambda\cap\{f\ne \hat\sigma \})\le \frac 1{6} \lambda^3.
\end{equation*}
We define $v\in SBD^2(Q_\lambda)$ by
\begin{equation}\label{eqdefvfrow}
 v(x):=w(x)-w(x+\frac12 e_1).
\end{equation}
Obviously $\calH^2(J_v \cap Q_\lambda)\le
\calH^2(J_w \cap Q_\lambda)+\calH^2(J_w \cap \hat Q_\lambda)
\le |De_{23}(u)|(Q_\lambda\cup \hat Q_\lambda)\le \eta\lambda^2$, and similarly
\begin{equation}
\int_{Q_\lambda} |e(v)|^2d\calL^3
\le 2 \int_{Q_\lambda\cup\hat Q_\lambda} |e(w)|^2d\calL^3
 \le 2 E_\eps[u, Q_\lambda\cup \hat Q_\lambda].
\end{equation}
We apply Lemma \ref{lemmaccfcube}\ref{lemmaccfcuberigid}  to $v$, and we obtain
a Borel set $\omega\subseteq Q_\lambda$ and an affine map $a:\R^3\to\R^3$ such that
\begin{equation*}
|\omega|\le \frac1{12} \lambda^3
\end{equation*}
and
\begin{equation}\label{eqwAffb}
 \int_{Q_\lambda\setminus\omega} |v-a|^2 d\calL^3\le c \lambda^2 E_\eps[u, Q_\lambda\cup \hat Q_\lambda]\le c \lambda^5 E.
\end{equation}
We then collect the three exceptional sets, and define
\begin{equation*}
 \tilde \omega:=\omega \cup (Q_\lambda\cap\{f\ne\sigma\})
 \cup ((\hat Q_\lambda\cap\{f\ne\hat \sigma\})-\frac12 e_1)
 ,\hskip2mm
 \text{ which obeys }\calL^3(\tilde\omega)\le\frac14\lambda^3.
\end{equation*}
From the definitions of $v$ and $w$ we obtain
\begin{equation*}\begin{split}
v_3(x)&=u_3(x)-u_3(x+\frac12 e_1),\\
v_2(x)&=u_2(x)-u_2(x+\frac12 e_1)+2(\hat\sigma-\sigma) x_3
\end{split}\end{equation*}
for almost all $x\in Q_\lambda\setminus\tilde\omega$; recalling the estimates for $u_2$ and $u_3$ in \eqref{eqchoicel2b} this leads to
\begin{equation*}
\begin{split}
   & \int_{Q_\lambda\setminus\tilde\omega}
 |v_3(x)+\frac12 \beta x_2 +\frac12 \bar t|^2
  d\calL^3 \le c   \lambda^3  E,\\ &\int_{Q_\lambda\setminus\tilde\omega}
  |v_2(x)+\frac12 \beta x_3  +\frac12 \bar s -2(\hat\sigma-\sigma) x_3|^2  d\calL^3 \le c   \lambda^3  E.
\end{split}
\end{equation*}
With \eqref{eqwAffb}, a triangular inequality and $\lambda\le1$ we see that the same two estimates hold
with $v$ replaced by the affine map $a$ from \eqref{eqwAffb}. We write it
in the form $a(x)=B+S\times x$, for some $B,S\in\R^3$, so that in particular
$a_2(x)=B_2+S_3x_1-S_1x_3$ and
$a_3(x)=B_3+S_1x_2-S_2x_1$, and obtain
\begin{equation*}
\begin{split}
    \int_{Q_\lambda\setminus\tilde\omega}
  |B_2+S_3x_1-S_1x_3+\frac12 \beta x_3  +\frac12 \bar s -2(\hat\sigma-\sigma) x_3|^2
  d\calL^3
  \le c   \lambda^3  E,\\
    \int_{Q_\lambda\setminus\tilde\omega}
 |B_3+S_1x_2-S_2x_1+\frac12 \beta x_2 +\frac12 \bar t|^2
  d\calL^3 \le c   \lambda^3  E.
\end{split}
\end{equation*}
By Lemma \ref{lemmacubelin}
we can estimate the coefficient of $x_3$ in the first integral, and the coefficient of $x_2$ in the second one. This leads to
\begin{equation*}
  \lambda^2 |-S_1+\frac12\beta -2(\hat\sigma-\sigma)|^2
+  \lambda^2 |S_1+\frac12\beta|^2
  \le c E
\end{equation*}
so that
\begin{equation}\label{eqalphasigsig}
  \lambda^2 |\beta -2(\hat\sigma-\sigma)|^2
   \le c  E
\end{equation}
for some $\hat\sigma,\sigma\in\{-1,1\}$.
As we assumed $|\beta|\le 3$, and
$\hat\sigma-\sigma\in\{-2,0,2\}$ we see that \eqref{eqalphasigsig} implies
\begin{equation}\label{eqenelastb}
  E \ge c \lambda^2\beta^2.
\end{equation}
Combining \eqref{eqesurfenb}
and \eqref{eqenelastb} we obtain
\begin{equation*}
 E\ge c\min\{ \frac{\eps}{\lambda}, \lambda^2\beta^2\} \hskip5mm
 \text{ for all }\lambda\in (0,\frac14].
\end{equation*}
If $\beta^2\ge\eps$ we choose $\lambda:=\frac14\eps^{1/3}|\beta|^{-2/3}$ and conclude
$E\ge c\eps^{2/3}|\beta|^{2/3}$.
If instead $\beta^2<\eps$ then
$\lambda=\frac14$ gives $E\ge c\beta^2$. This concludes the proof.
\end{proof}

\subsection{Lower bound for Neumann boundary data} \label{subsec:LB-6}

\begin{prop}\label{proplbneu}
There is $c>0$ such that for all $\eps\in (0,1]$ and all $\gamma\in\R$
\begin{equation*}
 c\min\{-\gamma^2, \frac{1}{c^2}\eps^{2/3}-|\gamma|\}\le \inf\{E_\eps[u]-\gamma M'(u): u\in W^{1,2}({\Omega;\R^3})\}.
\end{equation*}
\end{prop}
\begin{proof}
As usual, we fix $u\in W^{1,2}(\Omega;\R^3)$ and set $E:=E_\eps[u]$.
We can assume $E<\infty$ (as $M'(u)\in \R$ for all $u\in W^{1,2}(\Omega;\R^3)$).
In this proof, for clarity we give explicit names to many of the constants that appear in the various estimates.
We mostly denote by $c_X$ the (universal) constant entering the key estimate for quantity $X$.

As above, by Proposition~\ref{proplowerstructure} we obtain that there is $\beta\in \R$ with
\begin{equation}\label{eqbetae}
 |\beta|\le 1+c_\beta E^{1/2}+c_\beta E^{1/3}
\end{equation}
such that, for some $b_2$, $d_3$, $\bar s$ and $\bar t$ the estimates
\begin{equation}\label{eqestu2bbcc2}
\begin{split}
  \int_{\Omega} &
  \left[|u_2(x)-b_2(x_3)-\beta x_3 x_1 -\bar s x_1|^2
  \right.\\
  &+\left.
 |u_3(x)-d_3(x_2)-\beta x_1x_2 -\bar t x_1|^2
  \right]d\calL^3\le c E
\end{split}
\end{equation}
and
\begin{equation*}
 \int_{(-1,1)^2} \left| u_1(x_1,x_2,1)-u_1(x_1,x_2,-1)
 + 2\beta x_2+2\bar t\right| d\calL^2 \le c E^{1/2}
\end{equation*}
hold.
Recalling the definition of $M'(u)$ in \eqref{eqdefM1} and using
$\int_{(-1,1)^2} x_2^2 d\calL^2=\frac43$,
the last estimate implies
 \begin{equation}\label{eqMbetaE}
  \left|M'(u)+\frac83\beta\right| \le c_M E^{1/2}.
\end{equation}

We distinguish two cases. If $|\beta|\ge 2$, then \eqref{eqbetae} gives
\begin{equation*}
|\beta|\le 2 (|\beta|-1)\le 2c_\beta E^{1/2}+2c_\beta E^{1/3},
\end{equation*}
so that
$E\ge \min\{(2c_\beta)^{-2}, (2c_\beta)^{-3}\}$.
With
\eqref{eqMbetaE} we obtain
\begin{equation*}
 |M'(u)|\le \frac83 |\beta|+c_M E^{1/2}
 \le c'E^{1/2}.
\end{equation*}
Therefore
\begin{equation*}
  E-\gamma M'(u)\ge E - |\gamma| c' E^{1/2}\ge
  \min_{t\in\R} (t^2-\gamma c' t)=
  -\frac{(c')^2}{4}\gamma^2,
\end{equation*}
which concludes the proof for $|\beta|\ge2$.

Consider now the case $|\beta|< 2$. Then Proposition~\ref{proplocal} and \eqref{eqestu2bbcc2} give
\begin{equation*}%\label{eqEbetaloc}
  E\ge c_L\min\{\beta^2,\eps^{2/3}|\beta|^{2/3}\},
 \end{equation*}
where we can assume $c_L\le 1$.
Therefore, recalling that \eqref{eqMbetaE}
gives $|M'(u)|\le \frac83|\beta|+ c_M E^{1/2}$,
\begin{equation*}
\begin{split}
 E-\gamma M'(u)&\ge  \frac12 E-c_M|\gamma| E^{1/2}+
 \frac12 E-\frac 83|\gamma|\,|\beta|  \\
 &\ge \frac12 (E-2c_M|\gamma| E^{1/2})
 + \frac12  \left(c_L\min\{\beta^2,\eps^{2/3}|\beta|^{2/3}\}-\frac{16}3 |\gamma|\,|\beta|\right)\\
 &\ge \frac12 \min_{t\in\R } (t^2-2c_M\gamma t)
 + \frac12 \min_{0\le t\le 2} \left(c_L\min\{t^2,\eps^{2/3}t^{2/3}\}-\frac{16}3 |\gamma| t\right).
\end{split}
\end{equation*}
The first minimum is $-c_M^2\gamma^2$.
For the second one, we observe that $\min_t(c_L t^2-\frac{16}3 \gamma t)=-c_L'\gamma^2$ for some $c_L'>0$.
Consider now
$\min_{0\le t\le2} (c_L\eps^{2/3}t^{2/3}-\frac{16}3 |\gamma| t)$.
By concavity, this is attained either at $t=0$ or at $t=2$, and hence it equals
$\min\{0, c_L(2\eps)^{2/3}-\frac{32}3|\gamma|\}$. Collecting terms,
\begin{equation*}
\begin{split}
 E-\gamma M'(u)&\ge
 -\frac{c_M^2}{2}\gamma^2 +
\frac12 \min\{-c_L'\gamma^2,0, c_L(2\eps)^{2/3}-\frac{32}3|\gamma|\}.
\end{split}
\end{equation*}
As the first term is negative, we can drop 0 from the minimum. Further, the second term in the sum
controls the first, up to a factor $c_M^2/c_L'$. Therefore
\begin{equation*}
\begin{split}
 E-\gamma M'(u)&\ge
\frac12
(1+\frac{c_M^2}{c_L'})
\min\{-c_L'\gamma^2, c_L(2\eps)^{2/3}-\frac{32}3|\gamma|\} \\
&\ge  \min\{-c\gamma^2,\frac1{c}\eps^{2/3}-c|\gamma|\}.
\end{split}
\end{equation*}
As each of the three terms in the last expression is nonincreasing in $c$, we can take a unique constant (see also \eqref{eqneumonot}).
This concludes the proof.
\end{proof}

\section*{Acknowledgments}
This work was partially supported
by the Deutsche Forschungsgemeinschaft through project 211504053/SFB1060 and
project 441211072/SPP2256, and by the National Science Foundation through
grants {OISE-0967140, DMS-1311833, and DMS-2009746}.
The research of Oleksandr Misiats was supported by {the Simons Foundation through} Collaboration Grant for Mathematicians no. 854856.

\end{document}